\documentclass{article}
\usepackage{amsmath}
\usepackage{amssymb}
\usepackage{mathtools}
\usepackage{verbatim}
\usepackage{enumerate}
\usepackage{graphicx}
\usepackage{float}
\usepackage{overpic}
\usepackage{subfigure}
\usepackage{caption}
\usepackage{geometry}
\usepackage{algorithm}
\usepackage{algorithmic}
\usepackage{multirow}
\usepackage{tikz}
\usepackage{epstopdf}
\usepackage[section]{placeins}
\newcommand\dd{\,\mathrm{d}}
\newcommand\ii{\,\mathrm{i} \,}
\newcommand\bR{\mathbb{R}}

\newcommand\cI{\mathcal{I}}

\newcommand\pd[2]{\dfrac{\partial{#1}}{\partial{#2}}}

\newcommand{\nbr}{\tau}

\newcommand{\idxi}[2]{#1_{#2}}
\newcommand{\idxij}[3]{#1_{{#2}, {#3}}}

\begin{document}
\title{Sparsifying preconditioner for the time-harmonic Maxwell's equations}
\author{Fei Liu$^\sharp$ and Lexing Ying$^{\dagger\sharp}$\\
  $\dagger$ Department of Mathematics, Stanford University\\
  $\sharp$ Institute for Computational and Mathematical Engineering, Stanford University
}
\date{}
\maketitle

\begin{abstract}
  This paper presents the sparsifying preconditioner for the time-harmonic Maxwell's equations in
  the integral formulation. Following the work on sparsifying preconditioner for the
  Lippmann-Schwinger equation, this paper generalizes that approach from the scalar wave case to the
  vector case. The key idea is to construct a sparse approximation to the dense system by minimizing
  the non-local interactions in the integral equation, which allows for applying sparse linear
  solvers to reduce the computational cost. When combined with the standard GMRES solver, the number
  of preconditioned iterations remains small and essentially independent of the frequency. This
  suggests that, when the sparsifying preconditioner is adopted, solving the dense integral system
  can be done as efficiently as solving the sparse system from PDE discretization.
\end{abstract}



\section{Introduction}

This paper concerns the time-harmonic scattering problem for the Maxwell's equations with
inhomogeneous permittivity.  For electromagnetic scattering problems, the solution typically has a
highly oscillatory form especially when the wave number is large. Due to the Nyquist theorem, at
least a constant number of grid points is needed per wavelength to capture the pattern of the
solution. As a consequence, the number of unknowns could be huge in high frequency regime.


Common approaches to solve this problem involve discretizing the PDE with the finite difference or
the finite element methods. Well known schemes include the Yee grid
\cite{yee1966numerical,champagne2001fdfd} and the N\'ed\'elec curl-conforming finite element scheme
\cite{nedelec1980mixed}.  Exploiting the sparsity of the discretized system, the multifrontal method
or the nested dissection factorization
\cite{george1973nested,duff1983multifrontal,liu1992multifrontal} are generally applied in this
scenario, where the setup and the solve costs are $O(N^2)$ and $O(N^{4/3})$ respectively, which is a
huge advantage over the na\"ive Gaussian elimination. However, directly discretizing the PDE suffers
from the pollution effect \cite{babuska1997}. Higher order schemes could help to reduce the
pollution error, but the corresponding larger stencil supports will soon make the nested dissection
factorization no longer as effective.

Instead of solving the PDE form, one can solve the integral form of the equation. There are
  several advantages of doing that. First, the integral equation approach trades the dispersion
  error of the PDE approaches for the quadrature error, which can often be controlled by high order
  quadrature rules. Another advantage of solving the integral form is that, the boundary conditions
are dealt with more naturally, unlike for the PDE form where we need to seek for artificial
absorbing boundary conditions such as the PML \cite{berenger1994,johnson2008,chew1994}. Despite all
those advantages of the integral form, there is a notable drawback: the integral equation is dense,
thus sparse matrix techniques cannot be applied directly to save the computational cost.

Recently, the sparsifying preconditioner \cite{ying2015ls,ying2015sp} was developed to address the
efficiency of solving the integral form. It was originally designed for the scalar time-harmonic
wave equations such as the Lippmann-Schwinger equation and the time-harmonic Schr\"odinger
equation. The idea is to numerically transform the dense linear system into a sparse one by
minimizing the non-local interactions in the integral system. The solving process of the sparse
system then serves as a preconditioner for the dense integral system, where the iteration number is
essentially independent of the problem size. This paper demonstrates that this idea can be
generalized to the time-harmonic Maxwell's equations with suitable modifications. Consequently, the
integral form of the time-harmonic Maxwell's equations can be solved with a cost as cheap as solving
the PDE form up to a few preconditioned iterations.

Despite that our method solves the integral equation with the same order of cost for solving
  the PDE form, we rely on a key assumption that the medium needs to be smooth such that Nystr\"om
  discretization on a uniform Cartesian grid can be used to give reasonably accurate approximations
  and that FFT can be applied for the forward operator. For cases where the medium has sharp
  transitions, we refer the reader to \cite{beylkin2009fast}.

The rest of the paper is organized as follows. Section \ref{sec:problem} introduces the dense
integral equation of the time-harmonic Maxwell's equations. Section \ref{sec:sparsify} describes the
details of the sparsifying preconditioner for that equation. Numerical results in Section
\ref{sec:numerical} demonstrate the effectiveness of the preconditioner. Conclusions and future work
are given in Section \ref{sec:conclusion}.

\section{Problem formulation}
\label{sec:problem}

This section formalizes the problem we aim to solve. The goal is to solve the electromagnetic
scattering problem with inhomogeneous permittivity in isotropic media. Following the notations in
Chapter 9 of \cite{colton2012inverse}, let $\varepsilon = \varepsilon(x) > 0$ be the electric
permittivity, $\mu = \mu_0 > 0$ be the magnetic permeability and $\sigma = \sigma(x)$ be the
electric conductivity. We assume $\varepsilon(x) \equiv \varepsilon_0$ and $\sigma(x) \equiv 0$
outside some compact region $\Omega$. Note that the media is isotropic and $\varepsilon, \mu,
\sigma$ are all scalars. Under this setting, the time-harmonic Maxwell's equations can be written as
\begin{gather}
\begin{dcases}
\nabla\times E(x) - \ii k H(x) = 0,\\
\nabla\times H(x) + \ii k \, (1 - m(x)) E(x) = 0,
\end{dcases}
\label{eqn:Maxwell}
\end{gather}
where $k = \sqrt{\varepsilon_0\mu_0}\omega$ is a normalizing factor and $\omega$ is the angular
frequency. $m(x)$ is given by
\begin{gather*}
m(x) \coloneqq 1 - \dfrac{1}{\varepsilon_0} \left( \varepsilon(x) + \ii \dfrac{\sigma(x)}{\omega} \right),
\end{gather*}
where $m(x) \equiv 0$ outside $\Omega$.

Eliminating $H(x)$ in \eqref{eqn:Maxwell} gives the equation for $E(x)$
\begin{gather}
\nabla\times (\nabla\times E(x)) - k^2 (1 - m(x)) E(x) = 0.
\label{eqn:Maxwell_E}
\end{gather}
For the scattering problem, the electric field $E(x)$ consists of two parts: the incident field
$E^i(x)$ and the scattered field $E^s(x)$, where $E^i(x)$ is known and satisfies the time-harmonic
Maxwell's equations of the homogeneous background
\begin{gather*}
\nabla\times (\nabla\times E^i(x)) - k^2 E^i(x) = 0.
\end{gather*}
The goal is to solve the scattered field $E^s(x)$ by
\begin{gather*}
\nabla\times (\nabla\times (E^i(x) + E^s(x))) - k^2 (1 - m(x)) (E^i(x) + E^s(x)) = 0,
\end{gather*}
where $E^s(x)$ satisfies the Silver-M\"uller radiation condition \cite{muller2013foundations}
\begin{gather*}
\lim_{|x| \to \infty} (\nabla \times E^s(x))\times x - \ii k |x| E^s(x) = 0.
\end{gather*}
Following Chapter 9.2 of \cite{colton2012inverse}, an equivalent integral form of the equation is
given by
\begin{gather}
E(x) = E^i(x) - k^2 \int_{\bR^3} G(x - y) m(y) E(y) \dd y - \int_{\bR^3} \dfrac{1}{1 - m(y)} \nabla m(y) \cdot E(y) \, \nabla G(x - y) \dd y,
\label{eqn:Maxwell_int}
\end{gather}
where $G(x) \coloneqq \dfrac{e^{\ii k |x|}}{4\pi |x|}$ is the Helmholtz kernel. This paper aims to
solve \eqref{eqn:Maxwell_int} efficiently. Note that while \eqref{eqn:Maxwell_int} is posted on the whole 3D space, it only needs to be solved on $\Omega$ since $m(x)$ is compact supported in $\Omega$. We also note that an equation similar to \eqref{eqn:Maxwell_int}
is valid in 2D as well. The only difference is that $G(x)$ will be the 2D Helmholtz kernel. We shall
restrict our discussion below to the 3D case for simplicity and clarity. Nonetheless, the approach
works for 2D as well.

By rearranging the terms, we have the equation for $E^s(x)$
\begin{equation}
  E^s(x) + k^2 \int_{\Omega} G(x - y) m(y) E^s(y) \dd y + \int_{\Omega} \dfrac{1}{1 - m(y)} \nabla m(y) \cdot E^s (y) \, \nabla G(x - y) \dd y = g(x)
  \label{eqn:scattter}
\end{equation}
with $ g(x) \coloneqq - k^2 \int_{\Omega} G(x - y) m(y) E^i(y) \dd y - \int_{\Omega} \dfrac{1}{1 -
  m(y)} \nabla m(y) \cdot E^i (y) \, \nabla G(x - y) \dd y$. Now let $E^1(x), E^2(x), E^3(x)$ be the
three components of $E^s(x)$, and $g^1(x), g^2(x), g^3(x)$ be the components of $g(x)$, and
introduce
\begin{gather*}
p^1(x) \coloneqq \dfrac{1}{1 - m(x)}\pd{m}{x_1}(x), \quad
p^2(x) \coloneqq \dfrac{1}{1 - m(x)}\pd{m}{x_2}(x), \quad
p^3(x) \coloneqq \dfrac{1}{1 - m(x)}\pd{m}{x_3}(x), \quad
\\
G^1(x) \coloneqq \pd{G}{x_1}(x), \quad
G^2(x) \coloneqq \pd{G}{x_2}(x), \quad
G^3(x) \coloneqq \pd{G}{x_3}(x). \quad
\end{gather*}
With these notations, \eqref{eqn:scattter} can be rewritten as the following matrix form
\begin{gather*}
\begin{bmatrix}
E^1\\
E^2\\
E^3
\end{bmatrix}
+
k^2 G \ast \left( m
\begin{bmatrix}
E^1\\
E^2\\
E^3
\end{bmatrix}
\right)
+
\begin{bmatrix}
G^1 \ast\\
G^2 \ast\\
G^3 \ast
\end{bmatrix}
\left(
\begin{bmatrix}
p^1 & p^2 & p^3
\end{bmatrix}
\begin{bmatrix}
E^1\\
E^2\\
E^3
\end{bmatrix}
\right)
=
\begin{bmatrix}
g^1\\
g^2\\
g^3
\end{bmatrix}.
\end{gather*}

Without loss of generality, we assume that $\Omega = (0, 1)^3$ and discretize $\Omega$ with a
uniform Cartesian grid so that the convolutions can be evaluated efficiently by the FFT
\cite{cooley1965algorithm}. Let $n$ be the number of points per dimension and $h = 1/ (n+ 1)$ be the
step size. Denote
\begin{gather*}
\cI \coloneqq \{i = (i_1, i_2, i_3) : 1 \le i_1, i_2, i_3 \le n \}
\end{gather*}
as the discrete index set. To obtain the discrete equation, we use subscripts to denote the discrete
indices. For example, $\idxi{m}{i}$ stands for the value of $m(x)$ at $x = i h$ where $i = (i_1,
i_2, i_3) \in \cI$ is a multi-index. Then the discrete equations can be expressed as
\begin{gather}
\idxi{E^d}{i} + k^2 \sum_{j \in \cI} \idxij{G}{i}{j} \idxi{m}{j} \idxi{E^d}{j} + \sum_{j \in \cI} \idxij{G^d}{i}{j} (\idxi{p^1}{j} \idxi{E^1}{j} + \idxi{p^2}{j} \idxi{E^2}{j} + \idxi{p^3}{j} \idxi{E^3}{j}) = \idxi{g^d}{i}, \quad i \in \cI, \quad d = 1, 2, 3,
\label{eqn:discrete}
\end{gather}
where we slightly abuse the notation by using the same letter for the continuous and discrete
objects. To clarify, $\idxij{G}{i}{j}$ is the $(i,j)$-th entry of the corresponding convolution
(Toeplitz) matrix. For entries away from the diagonal, the values are given by
\begin{gather*}
\idxij{G}{i}{j} = h^3 G( ih - jh),
\end{gather*}
and for entries close to or on the diagonal where the Helmholtz kernel is singular, the values are
given by a fourth-order quadrature correction (see \cite{duan2009quadratures} for details). The same
notation is used for the partial derivative $G^d$ matrices. Higher order quadrature corrections can
also be used without modifying the following discussion.

Combining \eqref{eqn:discrete} for all $i\in\cI$ results in the discrete equation in matrix form
\begin{gather}
\begin{bmatrix}
E^1 \\
E^2 \\
E^3
\end{bmatrix}
+
\begin{bmatrix}
k^2 G & & & G^1 \\
& k^2 G & & G^2 \\
& & k^2 G & G^3
\end{bmatrix}
\begin{bmatrix}
m & & \\
& m & \\
& & m \\
p^1 & p^2 & p^3
\end{bmatrix}
\begin{bmatrix}
E^1 \\
E^2 \\
E^3
\end{bmatrix}
=
\begin{bmatrix}
g^1 \\
g^2 \\
g^3
\end{bmatrix}
\label{eqn:discrete_block}
\end{gather}
where $E^d$ and $g^d$ are discrete vectors. $m$ and $p^d$ should be interpreted as diagonal
matrices and $G$ and $G^d$ are convolution (Toeplitz) matrices.

\section{Sparsifying preconditioner}
\label{sec:sparsify}

To solve \eqref{eqn:discrete} efficiently, we adopt the idea of the sparsifying preconditioner
\cite{ying2015ls}. The key insight is that, as the integral equation comes from PDE formulation,
there exists some local stencil that can restrict any unknown to interact only to its nearby
neighbors. As a result, a sparse system can be formulated to approximate the dense one. Thereafter,
the process of solving the sparse system can be treated as a preconditioning step for the dense
system.

\subsection{Building the approximating sparse system}

For each $i \in \cI$, we denote $\nbr_i$ as its neighborhood
\begin{gather*}
\nbr_i \coloneqq \{ j \in \cI: \|j - i\|_\infty \le 1 \}.
\end{gather*}
Each index $i$ is involved with three unknowns $\idxi{E^1}{i}, \idxi{E^2}{i}$ and $\idxi{E^3}{i}$,
and the total number of unknowns is $3 n^3$. What we are going to do next is to construct three
equations for each $i$ where each equation only involves unknowns indexed by $j \in \nbr_i$, unlike
in \eqref{eqn:discrete} where each equation is dense. To start with, let us pull out the equations
indexed by $\nbr_i$ in \eqref{eqn:discrete} and rearrange them into the following form by splitting
the interactions into the local part (unknowns indexed by $\nbr_i$) and the non-local part (unknowns
indexed by $\nbr_i^c$):
\begin{equation}
\begin{gathered}
\idxi{E^d}{\nbr_i} + k^2 (\idxij{G}{\nbr_i}{\nbr_i}\idxi{(mE^d)}{\nbr_i} + \idxij{G}{\nbr_i}{\nbr_i^c}\idxi{(mE^d)}{\nbr_i^c}) \\
+ (\idxij{G^d}{\nbr_i}{\nbr_i} \idxi{(p^1 E^1 + p^2 E^2 + p^3 E^3)}{\nbr_i} + \idxij{G^d}{\nbr_i}{\nbr_i^c} \idxi{(p^1 E^1 + p^2 E^2 + p^3 E^3)}{\nbr_i^c}) = \idxi{g^d}{\nbr_i}, \quad d = 1, 2, 3.
\end{gathered}
 \label{eqn:nbr_i}
\end{equation}
We make the following explanations to clarify the notations in \eqref{eqn:nbr_i}:
\begin{itemize}
\item
$\nbr_i^c \coloneqq \cI \setminus \nbr_i$, which is the complement of $\nbr_i$ with respect to $\cI$.
\item
The single-subscript stands for the restriction of the corresponding vector to certain row
indices. For example, $\idxi{m}{\nbr_i}$ means the restriction of $m$ to $\nbr_i$.
\item
The double-subscript stands for the restriction of the corresponding matrix to certain row and
column indices. For example, $G_{\nbr_i, \nbr_i^c}$ is the sub-matrix of $G$ with row index set
$\nbr_i$ and column index set $\nbr_i^c$. The other notions for sub-matrix terms such as
$\idxij{G^d}{\nbr_i}{\nbr_i}$ should be interpreted similarly.
\end{itemize}
Equivalently we have the block matrix form
\begin{equation}
\begin{gathered}
\begin{bmatrix}
\idxi{E^1}{\nbr_i} \\
\idxi{E^2}{\nbr_i} \\
\idxi{E^3}{\nbr_i}
\end{bmatrix}
+
\begin{bmatrix}
k^2 \idxij{G}{\nbr_i}{\nbr_i} & & & \idxij{G^1}{\nbr_i}{\nbr_i} \\
& k^2 \idxij{G}{\nbr_i}{\nbr_i} & & \idxij{G^2}{\nbr_i}{\nbr_i} \\
& & k^2 \idxij{G}{\nbr_i}{\nbr_i} & \idxij{G^3}{\nbr_i}{\nbr_i} 
\end{bmatrix}
\begin{bmatrix}
\idxi{m}{\nbr_i} & & \\
& \idxi{m}{\nbr_i} & \\
& & \idxi{m}{\nbr_i} \\
\idxi{p^1}{\nbr_i} & \idxi{p^2}{\nbr_i} & \idxi{p^3}{\nbr_i}
\end{bmatrix}
\begin{bmatrix}
\idxi{E^1}{\nbr_i} \\
\idxi{E^2}{\nbr_i} \\
\idxi{E^3}{\nbr_i}
\end{bmatrix} \\
+
\begin{bmatrix}
k^2 \idxij{G}{\nbr_i}{\nbr_i^c} & & & \idxij{G^1}{\nbr_i}{\nbr_i^c} \\
& k^2 \idxij{G}{\nbr_i}{\nbr_i^c} & & \idxij{G^2}{\nbr_i}{\nbr_i^c} \\
& & k^2 \idxij{G}{\nbr_i}{\nbr_i^c} & \idxij{G^3}{\nbr_i}{\nbr_i^c} 
\end{bmatrix}
\begin{bmatrix}
\idxi{m}{\nbr_i^c} & & \\
& \idxi{m}{\nbr_i^c} & \\
& & \idxi{m}{\nbr_i^c} \\
\idxi{p^1}{\nbr_i^c} & \idxi{p^2}{\nbr_i^c} & \idxi{p^3}{\nbr_i^c}
\end{bmatrix}
\begin{bmatrix}
\idxi{E^1}{\nbr_i^c} \\
\idxi{E^2}{\nbr_i^c} \\
\idxi{E^3}{\nbr_i^c}
\end{bmatrix} 
=
\begin{bmatrix}
\idxi{g^1}{\nbr_i} \\
\idxi{g^2}{\nbr_i} \\
\idxi{g^3}{\nbr_i}
\end{bmatrix}.
 \label{eqn:nbr_i_block}
\end{gathered}
\end{equation}
The next step is to transform this equation set into three approximately sparse equations where the
non-local interactions can be neglected. Specifically, let $\alpha$ be a matrix of row size $3
\times |\nbr_i|$ and column size $3$. Multiplying $\alpha^T$ on both sides of
\eqref{eqn:nbr_i_block} gives
\begin{equation}
\begin{gathered}
\alpha^T
\begin{bmatrix}
\idxi{E^1}{\nbr_i} \\
\idxi{E^2}{\nbr_i} \\
\idxi{E^3}{\nbr_i}
\end{bmatrix} 
+
\alpha^T
\begin{bmatrix}
k^2 \idxij{G}{\nbr_i}{\nbr_i} & & & \idxij{G^1}{\nbr_i}{\nbr_i} \\
& k^2 \idxij{G}{\nbr_i}{\nbr_i} & & \idxij{G^2}{\nbr_i}{\nbr_i} \\
& & k^2 \idxij{G}{\nbr_i}{\nbr_i} & \idxij{G^3}{\nbr_i}{\nbr_i} 
\end{bmatrix}
\begin{bmatrix}
\idxi{m}{\nbr_i} & & \\
& \idxi{m}{\nbr_i} & \\
& & \idxi{m}{\nbr_i} \\
\idxi{p^1}{\nbr_i} & \idxi{p^2}{\nbr_i} & \idxi{p^3}{\nbr_i}
\end{bmatrix}
\begin{bmatrix}
\idxi{E^1}{\nbr_i} \\
\idxi{E^2}{\nbr_i} \\
\idxi{E^3}{\nbr_i}
\end{bmatrix} \\
+
\alpha^T
\begin{bmatrix}
k^2 \idxij{G}{\nbr_i}{\nbr_i^c} & & & \idxij{G^1}{\nbr_i}{\nbr_i^c} \\
& k^2 \idxij{G}{\nbr_i}{\nbr_i^c} & & \idxij{G^2}{\nbr_i}{\nbr_i^c} \\
& & k^2 \idxij{G}{\nbr_i}{\nbr_i^c} & \idxij{G^3}{\nbr_i}{\nbr_i^c} 
\end{bmatrix}
\begin{bmatrix}
\idxi{m}{\nbr_i^c} & & \\
& \idxi{m}{\nbr_i^c} & \\
& & \idxi{m}{\nbr_i^c} \\
\idxi{p^1}{\nbr_i^c} & \idxi{p^2}{\nbr_i^c} & \idxi{p^3}{\nbr_i^c}
\end{bmatrix}
\begin{bmatrix}
\idxi{E^1}{\nbr_i^c} \\
\idxi{E^2}{\nbr_i^c} \\
\idxi{E^3}{\nbr_i^c}
\end{bmatrix} 
=
\alpha^T
\begin{bmatrix}
\idxi{g^1}{\nbr_i} \\
\idxi{g^2}{\nbr_i} \\
\idxi{g^3}{\nbr_i}
\end{bmatrix}.
\label{eqn:nbr_i_block_alpha}
\end{gathered}
\end{equation}
If it is possible to find some non-trivial $\alpha$ such that
\begin{equation}
\alpha^T
\begin{bmatrix}
k^2 \idxij{G}{\nbr_i}{\nbr_i^c} & & & \idxij{G^1}{\nbr_i}{\nbr_i^c} \\
& k^2 \idxij{G}{\nbr_i}{\nbr_i^c} & & \idxij{G^2}{\nbr_i}{\nbr_i^c} \\
& & k^2 \idxij{G}{\nbr_i}{\nbr_i^c} & \idxij{G^3}{\nbr_i}{\nbr_i^c} 
\end{bmatrix}
\approx
0,
\label{eqn:approxzero}
\end{equation}
we can safely discard the terms involving the non-local interactions in
\eqref{eqn:nbr_i_block_alpha} and obtain
\begin{equation}
\begin{gathered}
\alpha^T
\begin{bmatrix}
\idxi{E^1}{\nbr_i} \\
\idxi{E^2}{\nbr_i} \\
\idxi{E^3}{\nbr_i}
\end{bmatrix} 
+
\beta^T
\begin{bmatrix}
\idxi{m}{\nbr_i} & & \\
& \idxi{m}{\nbr_i} & \\
& & \idxi{m}{\nbr_i} \\
\idxi{p^1}{\nbr_i} & \idxi{p^2}{\nbr_i} & \idxi{p^3}{\nbr_i}
\end{bmatrix}
\begin{bmatrix}
\idxi{E^1}{\nbr_i} \\
\idxi{E^2}{\nbr_i} \\
\idxi{E^3}{\nbr_i}
\end{bmatrix} 
\approx
\alpha^T
\begin{bmatrix}
\idxi{g^1}{\nbr_i} \\
\idxi{g^2}{\nbr_i} \\
\idxi{g^3}{\nbr_i}
\end{bmatrix},
\end{gathered}
\label{eqn:nbr_i_block_alpha_approx}
\end{equation}
where $\beta$ is computed by
\begin{gather}
\beta =
\begin{bmatrix}
k^2 \idxij{G}{\nbr_i}{\nbr_i} & & & \idxij{G^1}{\nbr_i}{\nbr_i} \\
& k^2 \idxij{G}{\nbr_i}{\nbr_i} & & \idxij{G^2}{\nbr_i}{\nbr_i} \\
& & k^2 \idxij{G}{\nbr_i}{\nbr_i} & \idxij{G^3}{\nbr_i}{\nbr_i} 
\end{bmatrix}^T
\alpha.
\label{eqn:alphabeta}
\end{gather}
Note that the unknowns involved in \eqref{eqn:nbr_i_block_alpha_approx} are all indexed by $\nbr_i$,
hence sparse and local. By repeating this for every index $i\in \cI$, one gets an approximately
sparse system that can be solved efficiently by sparse solvers, if treating the approximately equal
sign as strictly equal.

The key question is, does there exist such an $\alpha$ and how do we find it? To address this
question, let us examine the Helmholtz kernel $G(x)$ and its derivatives $G^d(x)$. The key
observation is that they satisfy the same Helmholtz equation at $x$ away from $0$. Specifically
\begin{align*}
  (-\Delta - k^2) G(x) = 0, \quad & x \ne 0, \\
  (-\Delta - k^2) G^d(x) = 0, \quad & x \ne 0, \quad d = 1, 2, 3.
\end{align*}
Since the row set $\nbr_i$ and the column set $\nbr_i^c$ for all matrices in \eqref{eqn:approxzero}
are naturally disjoint, there exists some local stencil $\gamma$, which is a column vector of size
$\nbr_i$ that can be thought of as a discretization of the operator $(-\Delta - k^2)$, such that
\begin{gather*}
\gamma^T
\begin{bmatrix}
\idxij{G}{\nbr_i}{\nbr_i^c} & \idxij{G^1}{\nbr_i}{\nbr_i^c} & \idxij{G^2}{\nbr_i}{\nbr_i^c} & \idxij{G^3}{\nbr_i}{\nbr_i^c}
\end{bmatrix} \approx 0.
\end{gather*}
In other words, the off-diagonal blocks of $G$ and $G^d$ can be simultaneously annihilated by
$\gamma$. Once $\gamma$ is ready, setting $\alpha$ as
\begin{gather}
\alpha
=
\begin{bmatrix}
\gamma & & \\
& \gamma & \\
& & \gamma
\end{bmatrix},
\label{eqn:alphagamma}
\end{gather}
gives rise to
\begin{gather*}
\alpha^T
\begin{bmatrix}
k^2 \idxij{G}{\nbr_i}{\nbr_i^c} & & & \idxij{G^1}{\nbr_i}{\nbr_i^c} \\
& k^2 \idxij{G}{\nbr_i}{\nbr_i^c} & & \idxij{G^2}{\nbr_i}{\nbr_i^c} \\
& & k^2 \idxij{G}{\nbr_i}{\nbr_i^c} & \idxij{G^3}{\nbr_i}{\nbr_i^c} 
\end{bmatrix} \\
=
\begin{bmatrix}
\gamma^T & & \\
& \gamma^T & \\
& & \gamma^T
\end{bmatrix}
\begin{bmatrix}
k^2 \idxij{G}{\nbr_i}{\nbr_i^c} & & & \idxij{G^1}{\nbr_i}{\nbr_i^c} \\
& k^2 \idxij{G}{\nbr_i}{\nbr_i^c} & & \idxij{G^2}{\nbr_i}{\nbr_i^c} \\
& & k^2 \idxij{G}{\nbr_i}{\nbr_i^c} & \idxij{G^3}{\nbr_i}{\nbr_i^c} 
\end{bmatrix}
\approx
0.
\end{gather*}
The above justifies the existence of such $\alpha$ for the interior index points. For the boundary
indices, $\alpha$ exists since one can construct local absorbing boundary conditions (ABCs) to
approximate the Silver-M\"uller radiation condition reasonably well.

In the actual implementation $\alpha$ is obtained in a more numerically way. More specifically, we
consider the following optimization problem
\begin{gather*}
\min_{\alpha^T \alpha = I} \|\alpha^T M\|_F
\end{gather*}
where
\begin{gather*}
M
\coloneqq
\begin{bmatrix}
k^2 \idxij{G}{\nbr_i}{\nbr_i^c} & & & \idxij{G^1}{\nbr_i}{\nbr_i^c} \\
& k^2 \idxij{G}{\nbr_i}{\nbr_i^c} & & \idxij{G^2}{\nbr_i}{\nbr_i^c} \\
& & k^2 \idxij{G}{\nbr_i}{\nbr_i^c} & \idxij{G^3}{\nbr_i}{\nbr_i^c} 
\end{bmatrix}.
\end{gather*}
The solution is given by the column-concatenation of the smallest three left singular vectors of
$M$, and it can be easily acquired by computing the SVD of $M$.

Once $\alpha$ is ready, we compute $\beta$ by \eqref{eqn:alphabeta} and form the three approximately
sparse equations in \eqref{eqn:nbr_i_block_alpha_approx}. Assembling all the equations for each $i
\in \cI$ and replacing ``$\approx$'' with ``$=$'' results in a sparse system, which can be solved
efficiently by the nested dissection algorithm. The solving process can be treated as a
preconditioner for the dense system \eqref{eqn:discrete_block}. As we shall see in Section
\ref{sec:numerical}, when combined with the standard GMRES solver, the preconditioner takes only a
few iterations to converge, where the rate is insensitive to the problem size.

\subsection{Exploiting translational invariance to compute the stencils}
This section is concerned with the efficient computation of the stencils $\alpha$. From the
discussion above, it seemed that we need to compute the SVD of $M$ for each index point $i$, which
could be costly. It turns out that these repetitive computations are not needed due to the
translational invariance of $G$ and $G^d$. To be specific, we categorize the index points $i \in
\cI$ into the following groups
\begin{itemize}
\item
  The interior index point $i = (i_1, i_2, i_3)$ where $2 \le i_1, i_2, i_3 \le n - 1$.
\item
  The face point $i$ where one of the $i_d$ is $1$ or $n$.
\item
  The edge point $i$ where two of the $i_d$ is $1$ or $n$.
\item
  The vertex point $i$ where all three $i_d$ is $1$ or $n$.
\end{itemize}
For the interior point $i$, we translate the neighborhood $\nbr_i$ to $\nbr$ as
\begin{gather*}
\nbr = \{j : -1 \le j_1, j_2, j_3 \le 1\},
\end{gather*}
i.e., the neighborhood of the original point, and we set $\nbr^c$ as
\begin{gather*}
\nbr^c = \{j : -n + 2 \le j_1, j_2, j_3 \le n - 2 \} \setminus \nbr,
\end{gather*}
then we compute $\alpha$ and $\beta$ from the matrices
\begin{gather*}
\begin{bmatrix}
k^2 \idxij{G}{\nbr}{\nbr} & & & \idxij{G^1}{\nbr}{\nbr} \\
& k^2 \idxij{G}{\nbr}{\nbr} & & \idxij{G^2}{\nbr}{\nbr} \\
& & k^2 \idxij{G}{\nbr}{\nbr} & \idxij{G^3}{\nbr}{\nbr} 
\end{bmatrix}
\quad \text{and} \quad
\begin{bmatrix}
k^2 \idxij{G}{\nbr}{\nbr^c} & & & \idxij{G^1}{\nbr}{\nbr^c} \\
& k^2 \idxij{G}{\nbr}{\nbr^c} & & \idxij{G^2}{\nbr}{\nbr^c} \\
& & k^2 \idxij{G}{\nbr}{\nbr^c} & \idxij{G^3}{\nbr}{\nbr^c} 
\end{bmatrix}.
\end{gather*}
Due to the translational invariance property of the convolution matrices, the stencils $\alpha$ and
$\beta$ acquired here work for all the interior index points. Note that the complement $\nbr^c$ is
taken with respect to a larger index set so that each translated copy of $\nbr_i^c$ is covered.

For the face points, let us take $i = (1, i_2, i_3)$ as an example where $ 2 \le i_2, i_3 \le n -
1$. In this case, we should set
\begin{gather*}
  \nbr = \{j : 1 \le j_1 \le 2 \text{ and } -1 \le j_2, j_3 \le 1 \},\\
  \nbr^c = \{j : 1 \le j_1 \le n \text{ and } -n + 2 \le j_2, j_3 \le n - 2 \} \setminus \nbr,
\end{gather*}
and the rest of the procedure is the same as for the interior points.

For the edge points and the vertex points, the process above can be generalized without
difficulty. Take $i = (1, 1, i_3)$ where $2 \le i_3 \le n - 1$ for the edge point
example. Correspondingly, we set
\begin{gather*}
\nbr = \{j : 1 \le j_1, j_2 \le 2 \text{ and } -1 \le j_3 \le 1 \},\\
\nbr^c = \{j : 1 \le j_1, j_2 \le n \text{ and } -n + 2 \le j_3 \le n - 2 \} \setminus \nbr.
\end{gather*}

The last example is for the vertex point $i = (1,1,1)$ where we have
\begin{gather*}
\nbr = \{j : 1 \le j_1, j_2, j_3 \le 2 \},\\
\nbr^c = \{j : 1 \le j_1, j_2, j_3 \le n\} \setminus \nbr.
\end{gather*}

\subsection{Complexity analysis}
Let $N = 3 n^3$ be the number of unknowns. From the previous discussions, we see that computing the
stencils $\alpha$ and $\beta$ for all the index groups needs $O(N)$ time and $O(N)$ space in
total. Once we have the stencils, the sparse system can be formed and the nested dissection
algorithm can be applied. For this stage, the setup cost is $O(N^2)$ time and $O(N^{4/3})$ space,
and the application time cost is $O(N^{4/3})$ in 3D. The forward operator of the dense system can be
evaluated fast by the FFT with $O(N \log N)$ cost, dominated by the nested dissection
algorithm. Thus the overall costs are: $O(N^2)$ time and $O(N^{4/3})$ space for the preconditioner
setup, and $O(N^{4/3})$ time per preconditioner application.

In the 2D case, the setup cost is $O(N^{3/2})$ time and $O(N \log N)$ space, and the application
time cost is $O(N \log N)$. As shown by the numerical results in Section \ref{sec:numerical}, the
preconditioner converges in only a few iterations, essentially independent of the problem size and
the frequency. This implies that, by applying the sparsifying preconditioner, solving the dense
integral system is comparable to the cost of solving the sparse system.

\section{Numerical results}
\label{sec:numerical}

This section presents the numerical results. The algorithm is implemented in MATLAB and the tests
are performed on a server with four Intel Xeon E7-4830-v3 CPUs. The preconditioner is combined with
the standard GMRES solver. The relative residual is $10^{-6}$ and the restart value is 20. The step
size $h$ is chosen such that there are six points per background wavelength. Numerical examples are
presented in both 2D and 3D.

\paragraph{2D problems.} Three examples are considered, where the $m(x)$ is
\begin{enumerate}
\item
  \label{m:2D1}
  a converging Gaussian lens,
\item
  \label{m:2D2}
  a square obstacle with smooth boundary,
\item
  \label{m:2D3}
  a random perturbation of the square obstacle,
\end{enumerate}
respectively. The incident field $E^i(x)$ is a plane wave
\begin{gather*}
  E^i(x) = 
  \begin{bmatrix}
    0 \\
    e^{\ii k x_1}
\end{bmatrix}.
\end{gather*}

\begin{table}[!h]
  \centering
  \includegraphics[width=0.42\textwidth]{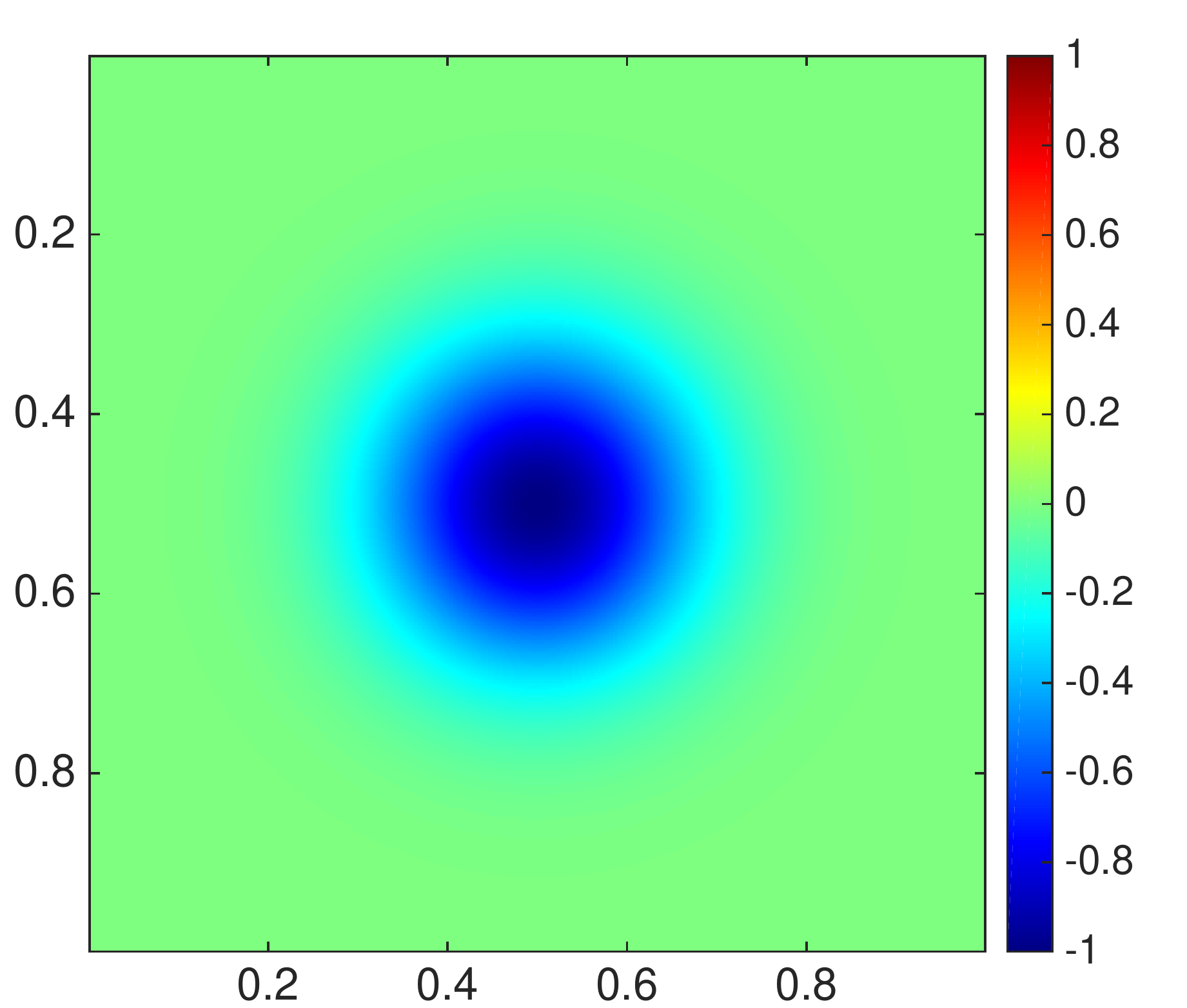}
  \includegraphics[width=0.42\textwidth]{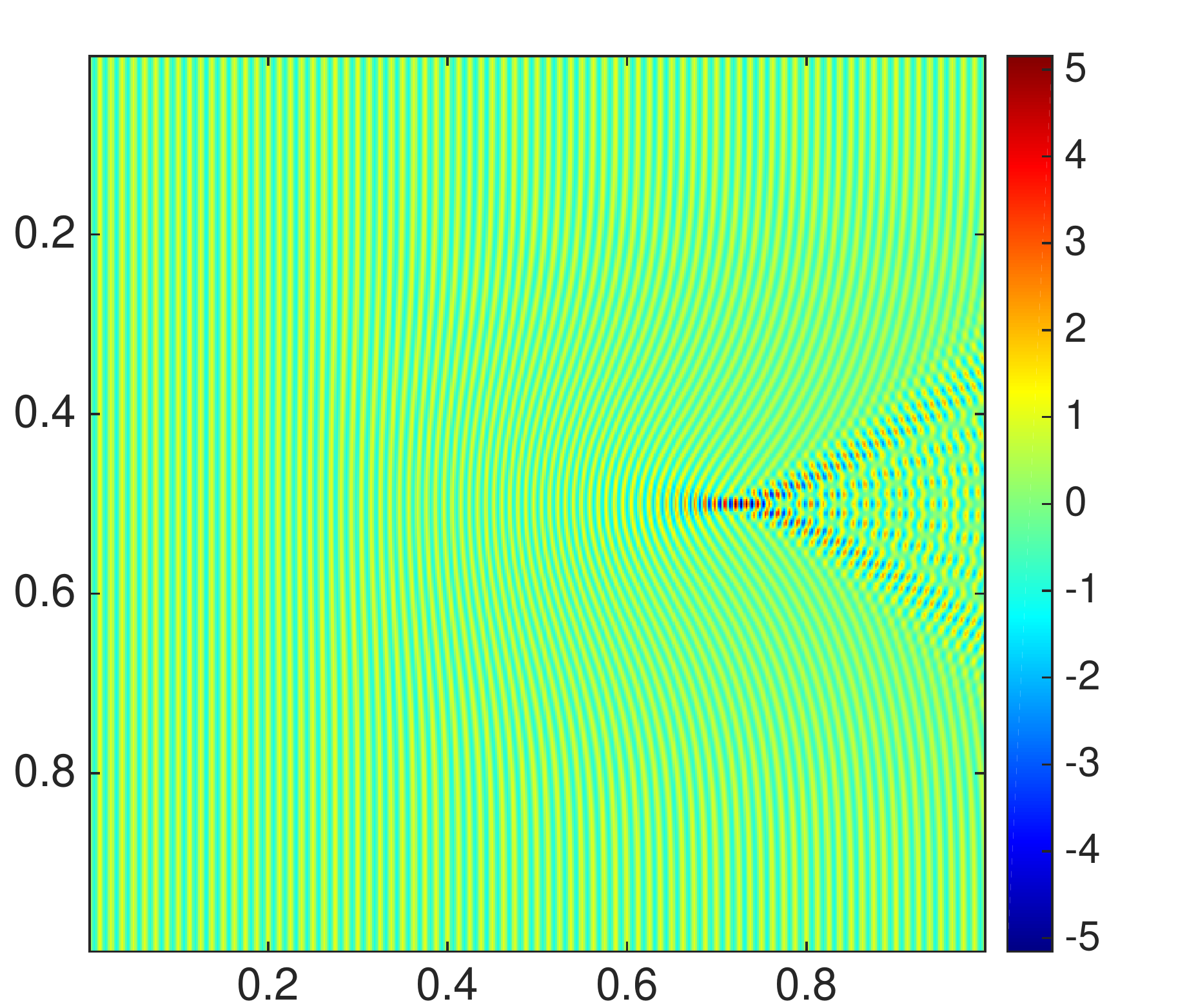}
  \begin{tabular}{cc | cc | cc}
    \hline
    \hline
    $k/(2\pi)$ & $N$ & $T_{\text{setup}}$ & $T_{\text{apply}}$ & $N_{\text{iter}}$ & $T_{\text{solve}}$ \\
    \hline
    $20$ & $2 \times 119^2$ & 1.25e+00 & 5.58e-02 & 6 & 9.13e-01 \\ 
    $40$ & $2 \times 239^2$ & 5.52e+00 & 3.27e-01 & 6 & 2.51e+00 \\ 
    $80$ & $2 \times 479^2$ & 2.17e+01 & 1.40e+00 & 6 & 8.93e+00 \\ 
    $160$ & $2 \times 959^2$ & 9.87e+01 & 5.03e+00 & 6 & 4.13e+01 \\ 
    \hline
    \hline
  \end{tabular}
  \caption{Results for example \eqref{m:2D1} in 2D. Top left: the inhomogeneity
    $m(x)$. Top right: the second component of the total field $E^i(x) + E^s(x)$ for $k/(2\pi) = 80$. Bottom: the
    numerical results.}
  \label{tab:2D1}
\end{table}

\begin{table}[!h]
  \centering
  \includegraphics[width=0.42\textwidth]{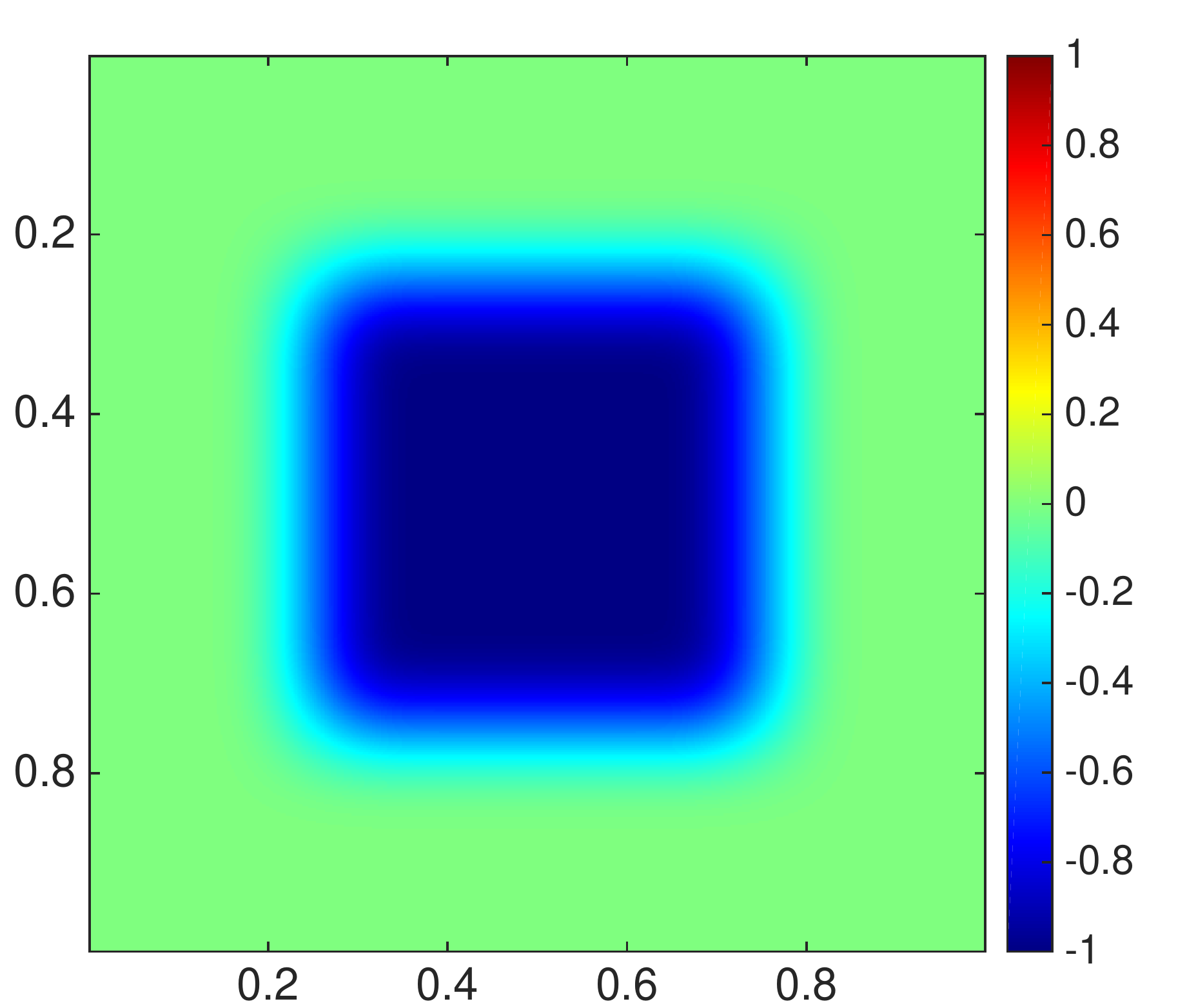}
  \includegraphics[width=0.42\textwidth]{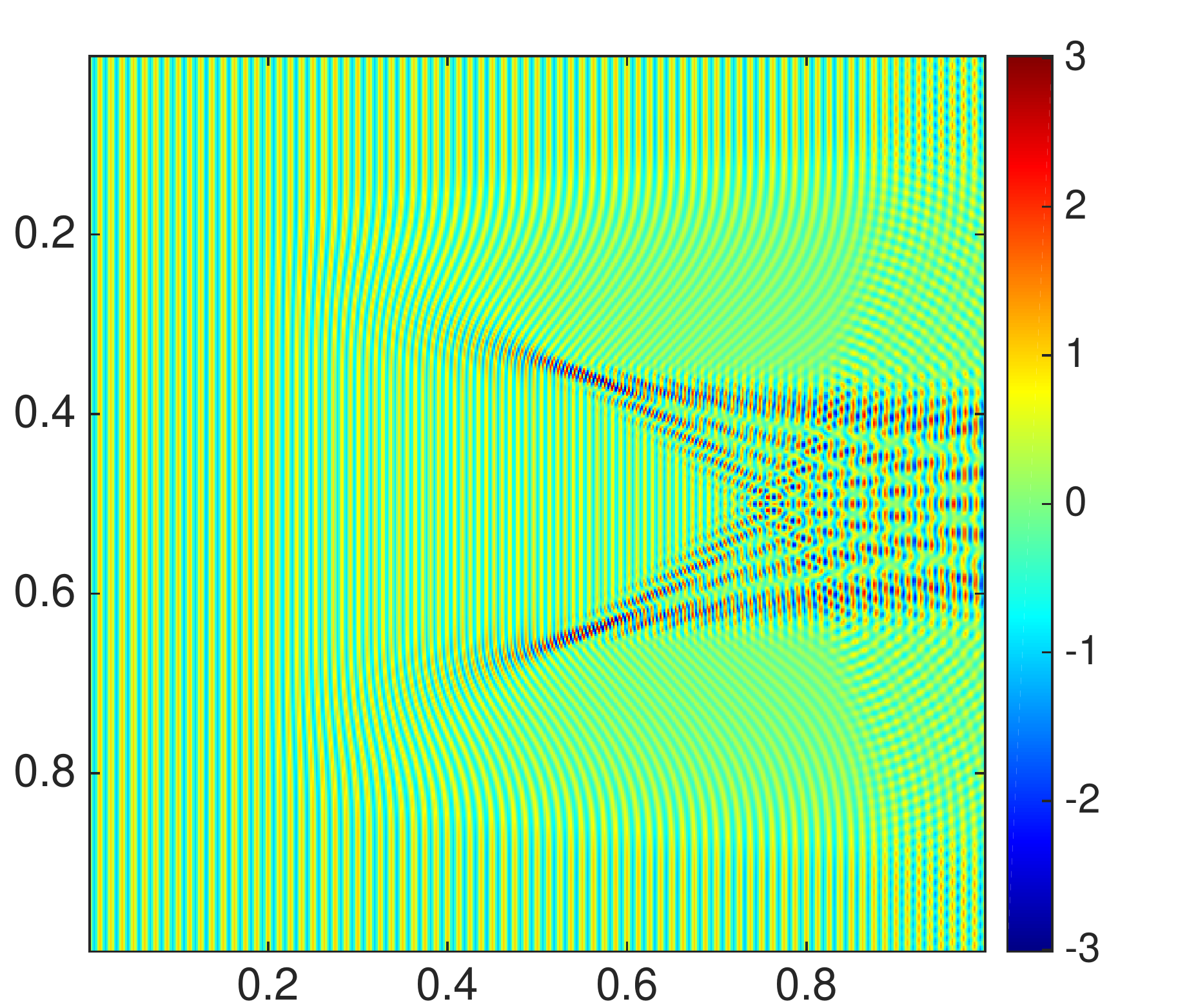}
  \begin{tabular}{cc | cc | cc}
    \hline
    \hline
    $k/(2\pi)$ & $N$ & $T_{\text{setup}}$ & $T_{\text{apply}}$ & $N_{\text{iter}}$ & $T_{\text{solve}}$ \\
    \hline
    $20$ & $2 \times 119^2$ & 1.28e+00 & 5.60e-02 & 6 & 8.41e-01 \\ 
    $40$ & $2 \times 239^2$ & 5.31e+00 & 2.46e-01 & 7 & 2.31e+00 \\ 
    $80$ & $2 \times 479^2$ & 2.24e+01 & 1.37e+00 & 8 & 1.32e+01 \\ 
    $160$ & $2 \times 959^2$ & 9.76e+01 & 5.85e+00 & 9 & 5.89e+01 \\ 
    \hline
    \hline
  \end{tabular}
  \caption{Results for example \eqref{m:2D2} in 2D. Top left: the inhomogeneity $m(x)$. Top right:
    the second component of the total field $E^i(x) + E^s(x)$ for $k/(2\pi) = 80$. Bottom: the
    numerical results.}
  \label{tab:2D2}
\end{table}

\begin{table}[!h]
  \centering
  \includegraphics[width=0.42\textwidth]{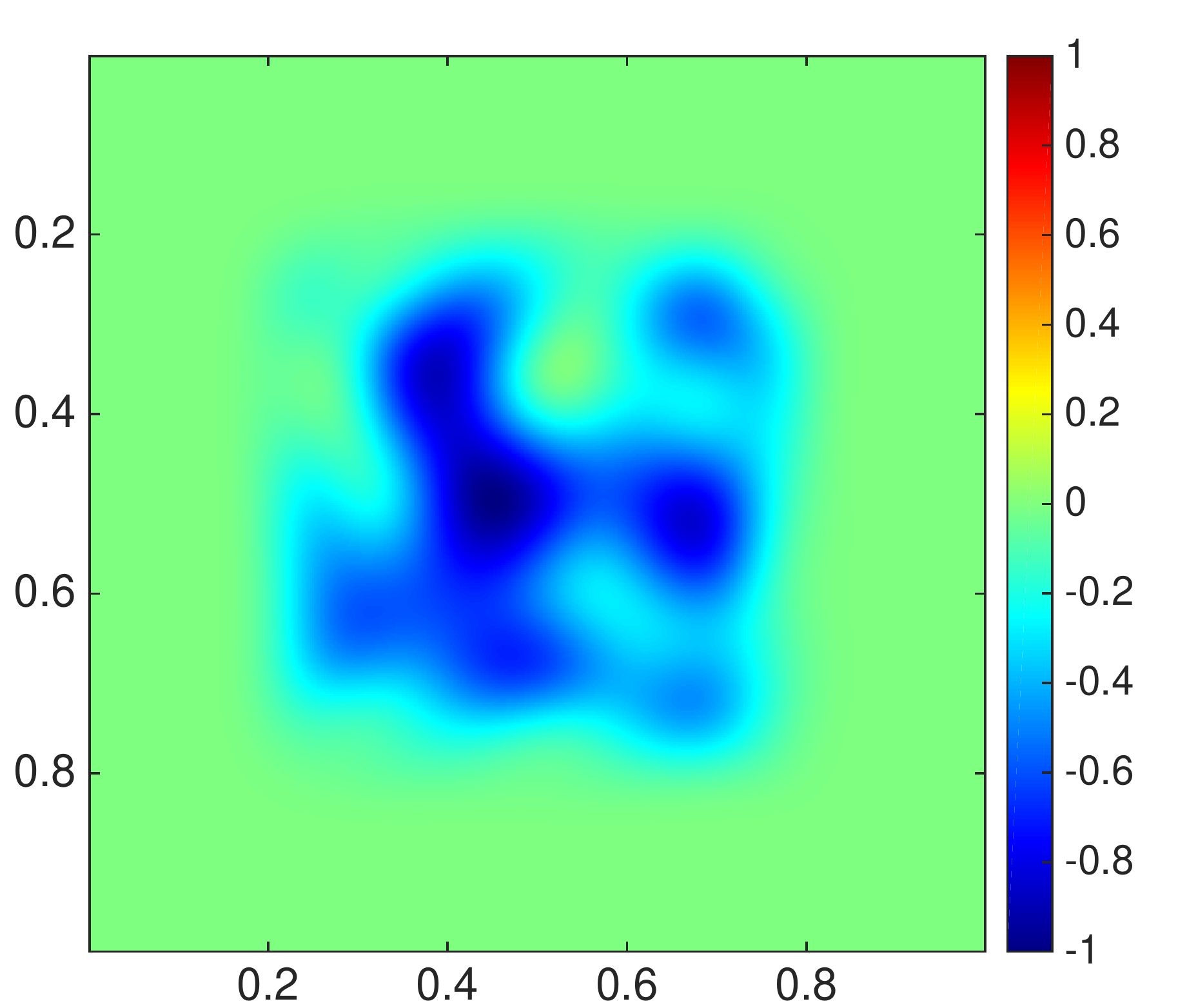}
  \includegraphics[width=0.42\textwidth]{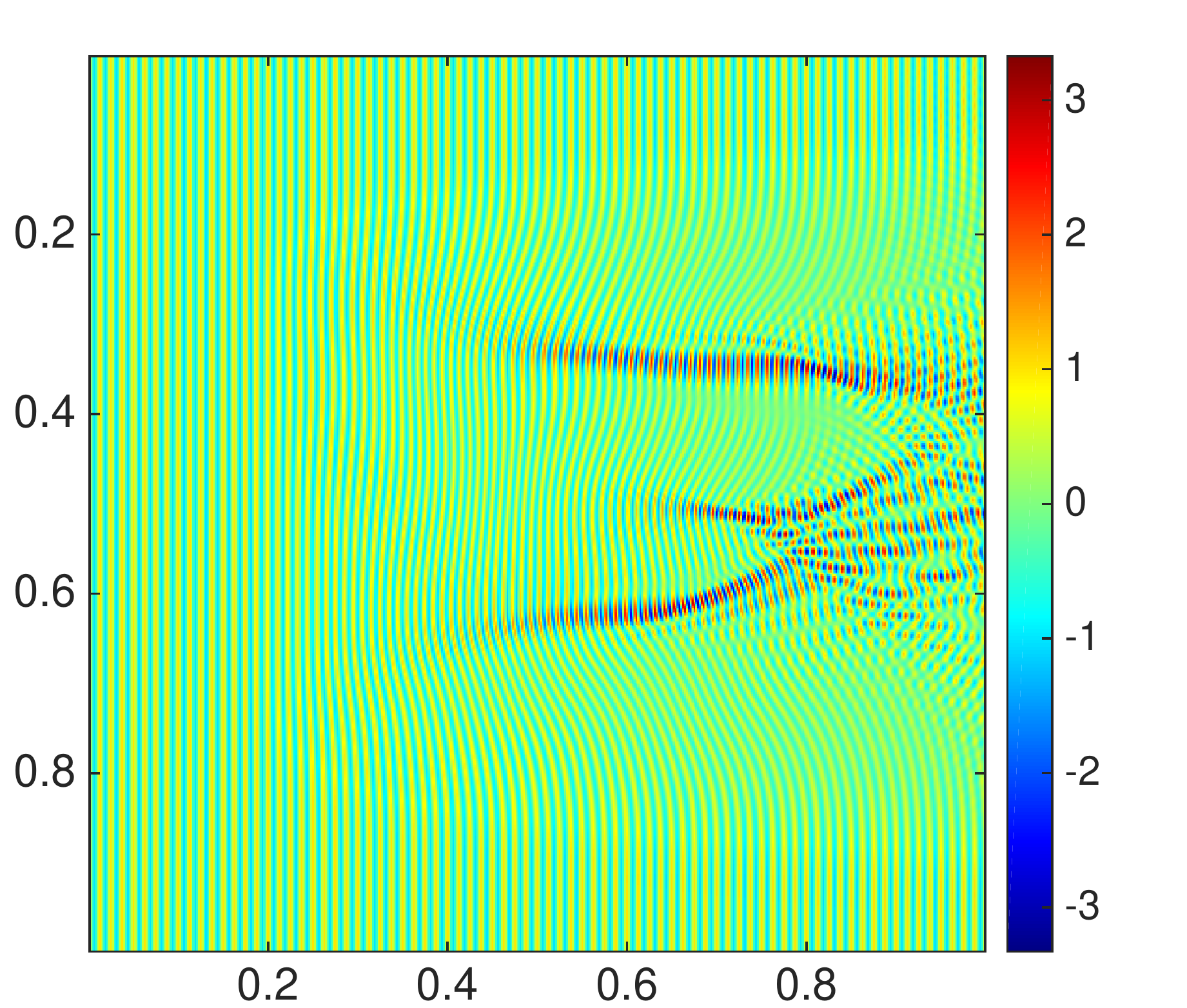}
  \begin{tabular}{cc | cc | cc}
    \hline
    \hline
    $k/(2\pi)$ & $N$ & $T_{\text{setup}}$ & $T_{\text{apply}}$ & $N_{\text{iter}}$ & $T_{\text{solve}}$ \\
    \hline
    $20$ & $2 \times 119^2$ & 1.38e+00 & 7.12e-02 & 6 & 1.19e+00 \\ 
    $40$ & $2 \times 239^2$ & 5.62e+00 & 2.68e-01 & 6 & 2.26e+00 \\ 
    $80$ & $2 \times 479^2$ & 2.30e+01 & 1.53e+00 & 6 & 9.83e+00 \\ 
    $160$ & $2 \times 959^2$ & 9.98e+01 & 6.36e+00 & 7 & 4.71e+01 \\ 
    \hline
    \hline
  \end{tabular}
  \caption{Results for example \eqref{m:2D3} in 2D. Top left: the inhomogeneity
    $m(x)$. Top right: the second component of the total field $E^i(x) + E^s(x)$ for $k/(2\pi) = 80$. Bottom: the
    numerical results.}
\label{tab:2D3}
\end{table}

The results of these three examples are given in Tables \ref{tab:2D1}, \ref{tab:2D2} and
\ref{tab:2D3}, respectively. The notations used in the tables are listed as follows:
\begin{itemize}
\item
$k / (2\pi)$ is the background wave number.
\item
$N$ is the number of unknowns.
\item
$T_{\text{setup}}$ is the setup cost of the preconditioner in seconds.
\item
$T_{\text{apply}}$ is the application cost of the preconditioner in seconds.
\item
$N_{\text{iter}}$ is the iteration number.
\item
$T_{\text{solve}}$ is the solve cost of the preconditioner in seconds.
\end{itemize}

\paragraph{3D problems.} Three examples are considered again, where the $m(x)$ is
\begin{enumerate}
\item
  \label{m:3D1}
  a converging Gaussian lens,
\item
  \label{m:3D2}
  a cube obstacle with smooth boundary,
\item
  \label{m:3D3}
  a random perturbation of the cube obstacle,
\end{enumerate}
respectively. The incident field $E^i(x)$ is
\begin{gather*}
E^i(x) = 
\begin{bmatrix}
0 \\
0 \\
e^{\ii k x_1}
\end{bmatrix}.
\end{gather*}
The results are given in Tables \ref{tab:3D1}, \ref{tab:3D2} and \ref{tab:3D3}, respectively.

\begin{table}[!ht]
\centering
\includegraphics[width=0.42\textwidth]{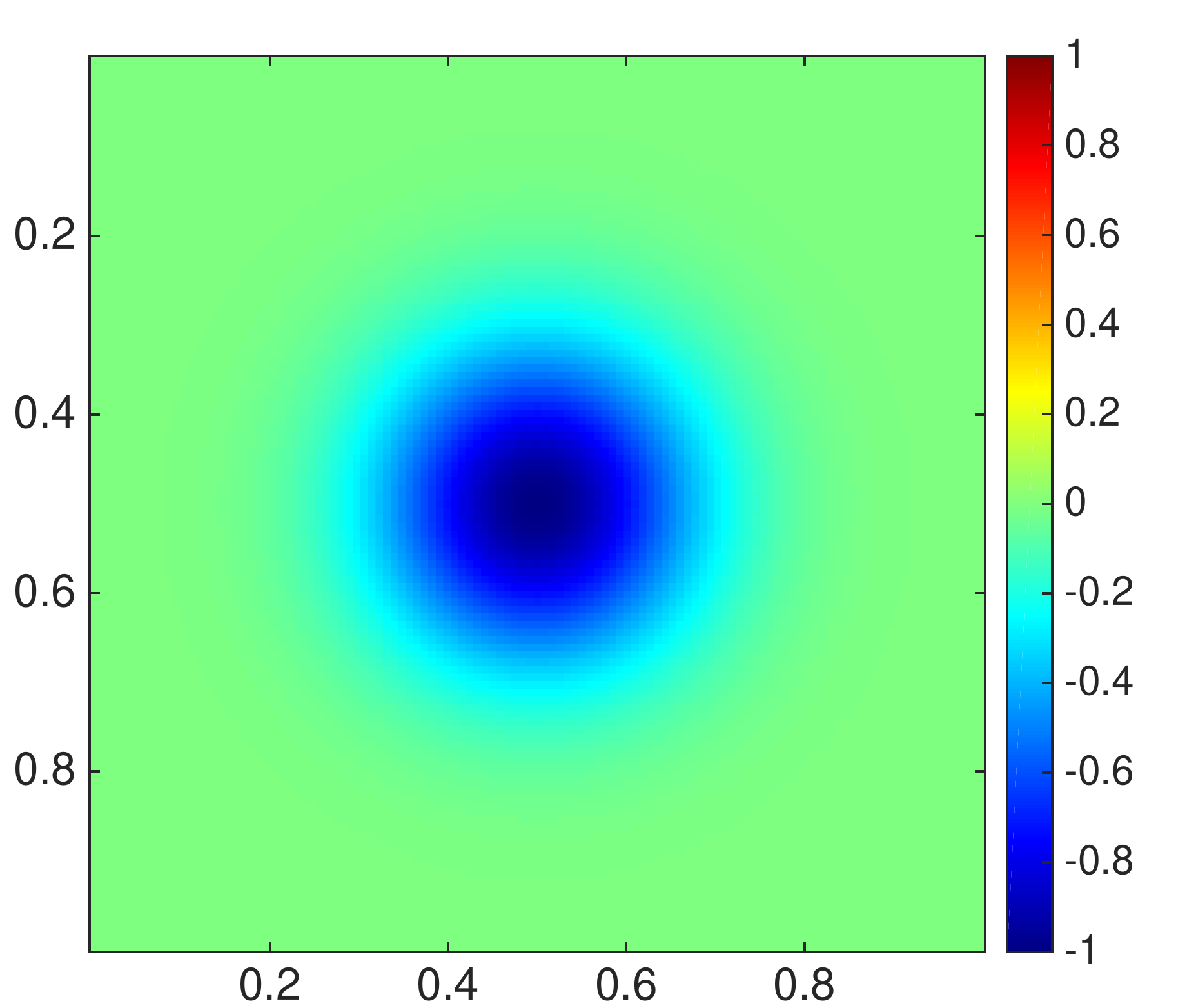}
\includegraphics[width=0.42\textwidth]{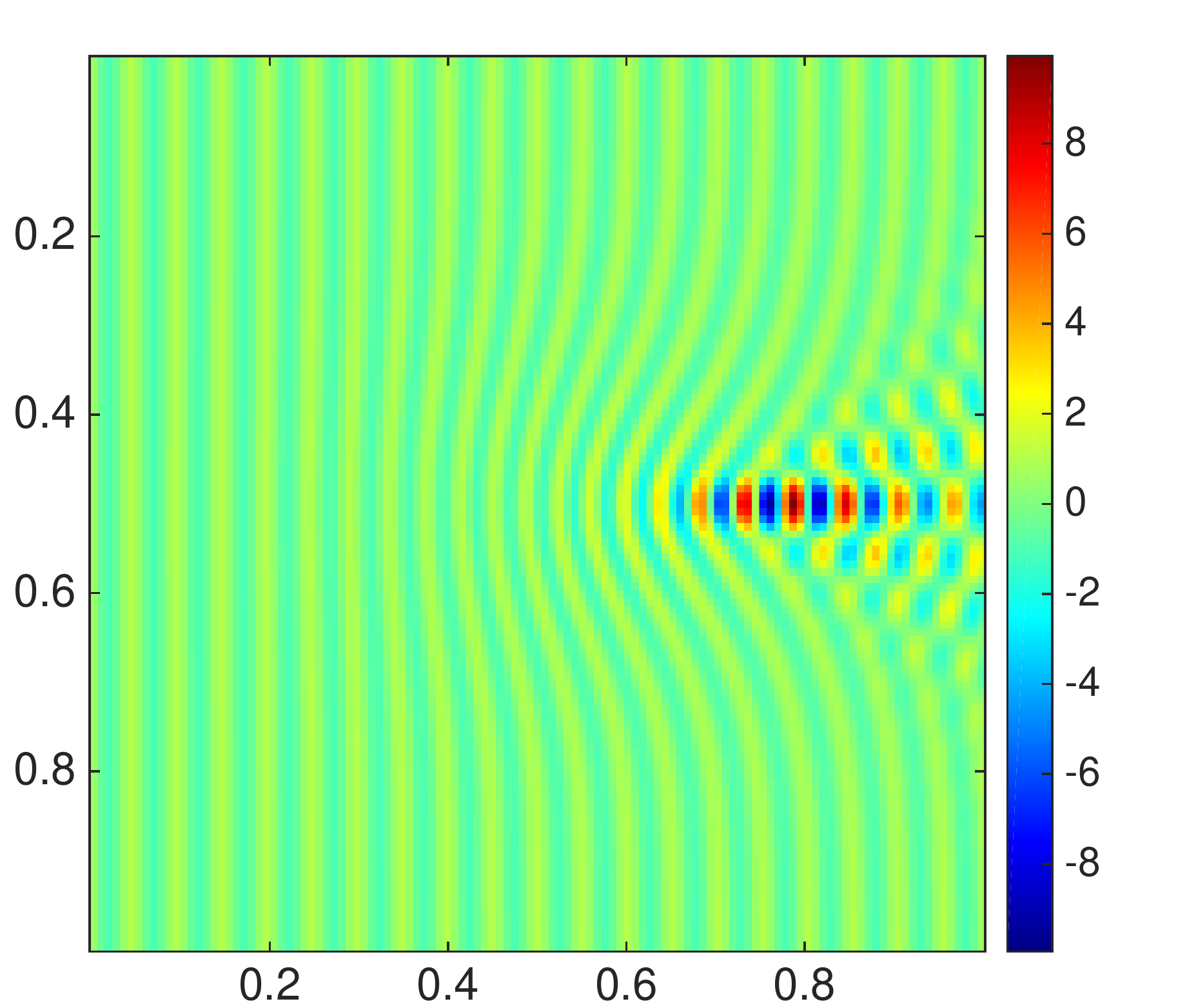}
\begin{tabular}{cc | cc | cc}
\hline
\hline
$k/(2\pi)$ & $N$ & $T_{\text{setup}}$ & $T_{\text{apply}}$ & $N_{\text{iter}}$ & $T_{\text{solve}}$ \\
\hline
$5$ & $3 \times 29^3$ & 2.07e+01 & 6.39e-01 & 6 & 5.10e+00 \\ 
$10$ & $3 \times 59^3$ & 4.84e+02 & 6.56e+00 & 6 & 4.67e+01 \\ 
$20$ & $3 \times 119^3$ & 1.18e+04 & 7.81e+01 & 6 & 5.31e+02 \\ 
\hline
\hline
\end{tabular}
\caption{Results for example \eqref{m:3D1} in 3D. Top left: the inhomogeneity $m(x)$ in
  cross section view at $x_3 = 0.5$. Top right: the third component of the total field $E^i(x) +
  E^s(x)$ at $x_3 = 0.5$ for $k/(2\pi) = 20$. Bottom: the numerical results.}
\label{tab:3D1}
\end{table}

\begin{table}[!ht]
\centering
\includegraphics[width=0.42\textwidth]{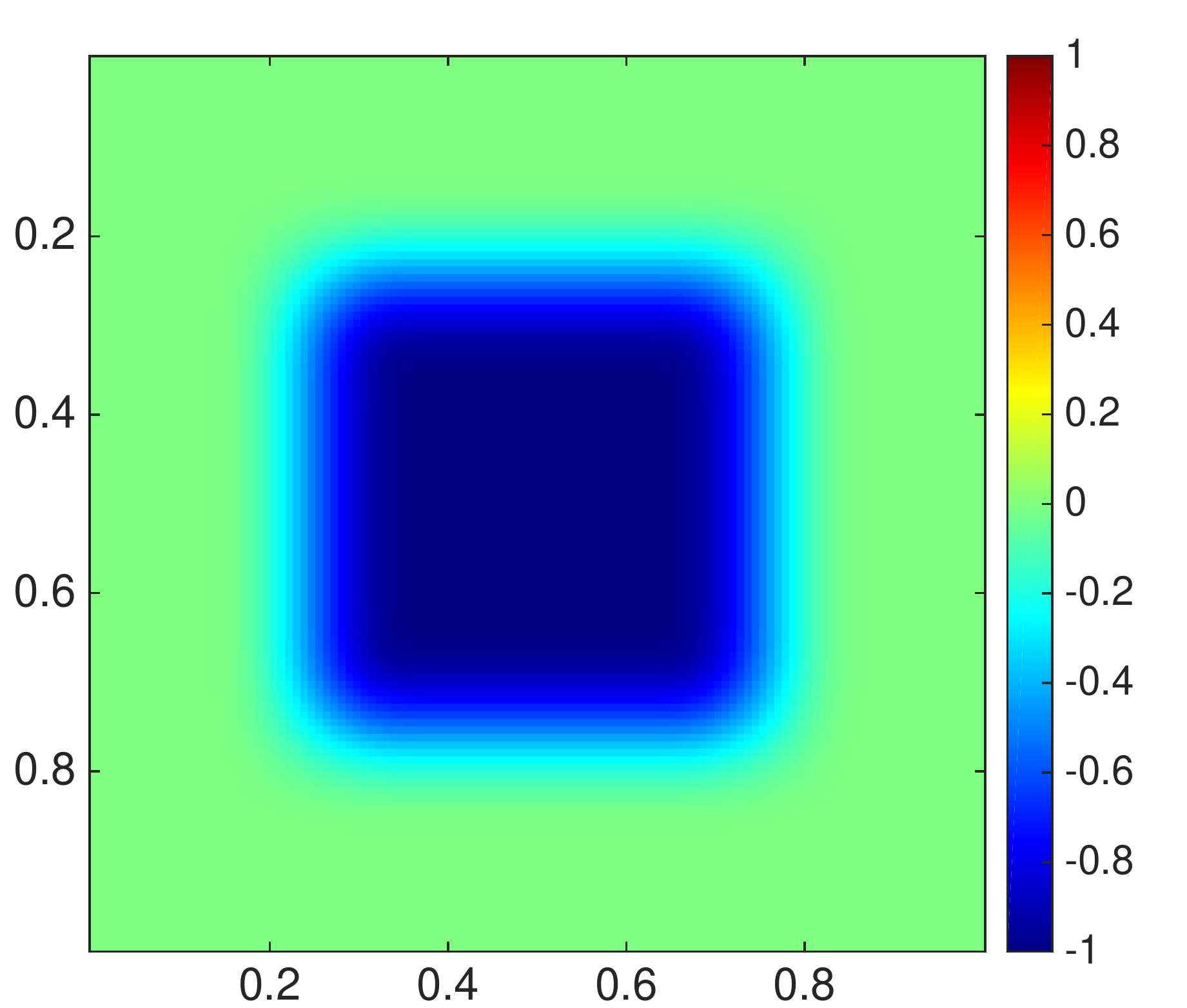}
\includegraphics[width=0.42\textwidth]{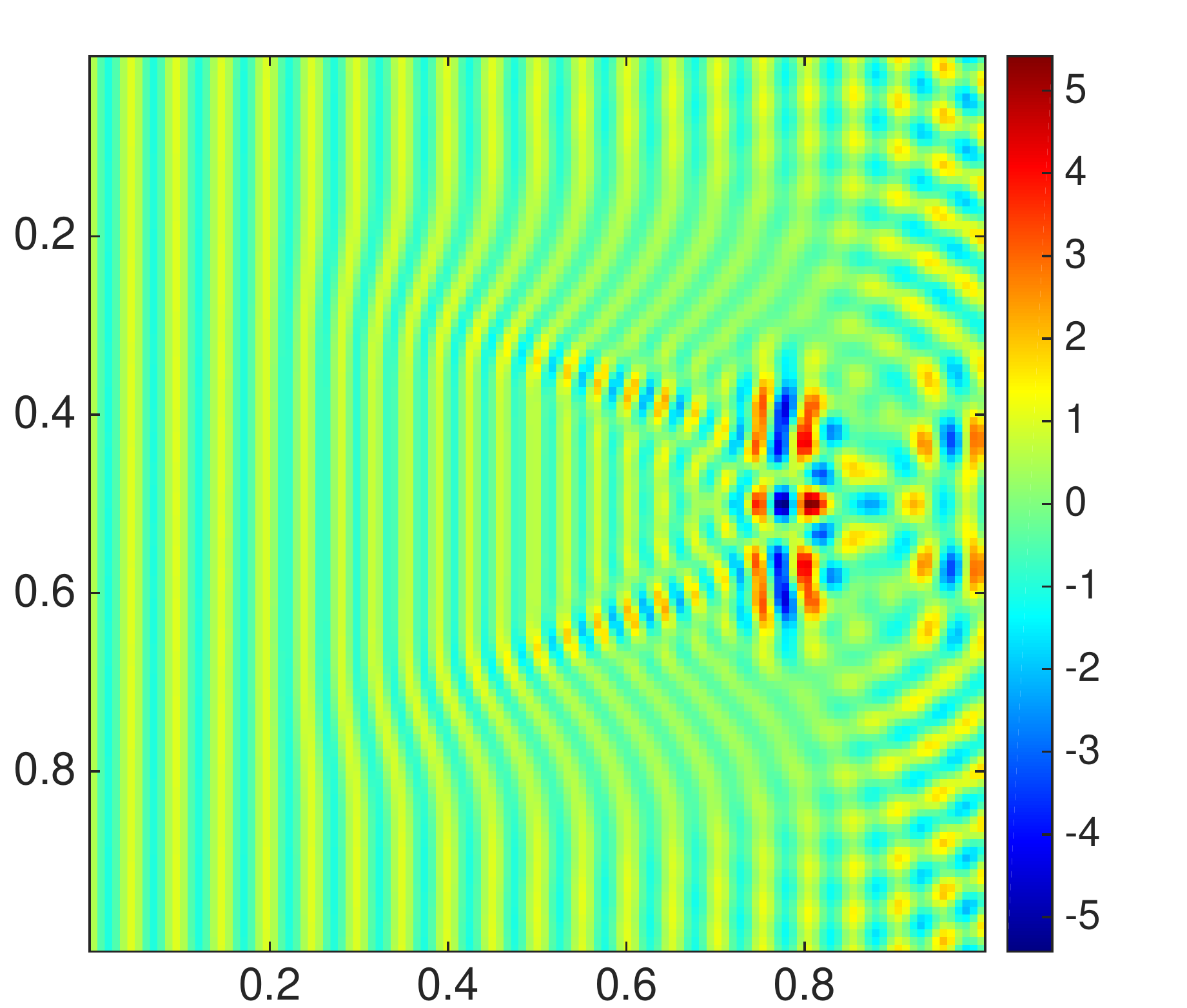}
\begin{tabular}{cc | cc | cc}
\hline
\hline
$k/(2\pi)$ & $N$ & $T_{\text{setup}}$ & $T_{\text{apply}}$ & $N_{\text{iter}}$ & $T_{\text{solve}}$ \\
\hline
$5$ & $3 \times 29^3$ & 2.17e+01 & 5.48e-01 & 7 & 5.18e+00 \\ 
$10$ & $3 \times 59^3$ & 4.83e+02 & 7.49e+00 & 7 & 5.52e+01 \\ 
$20$ & $3 \times 119^3$ & 1.18e+04 & 8.06e+01 & 7 & 6.20e+02 \\ 
\hline
\hline
\end{tabular}
\caption{Results for example \eqref{m:3D2} in 3D. Top left: the inhomogeneity $m(x)$ in
  cross section view at $x_3 = 0.5$. Top right: the third component of the total field $E^i(x) +
  E^s(x)$ at $x_3 = 0.5$ for $k/(2\pi) = 20$. Bottom: the numerical results.}
\label{tab:3D2}
\end{table}

\begin{table}[!ht]
\centering
\includegraphics[width=0.42\textwidth]{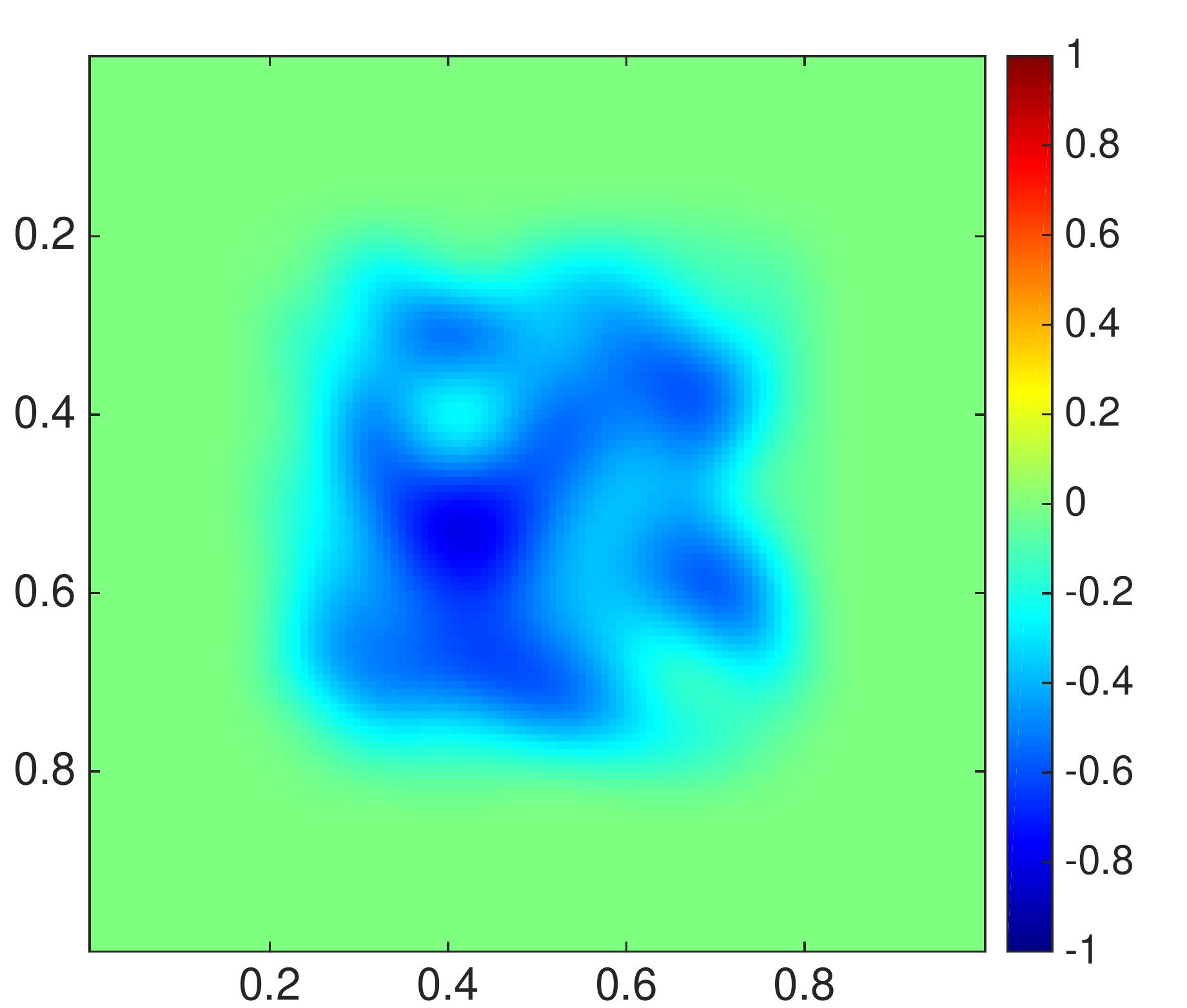}
\includegraphics[width=0.42\textwidth]{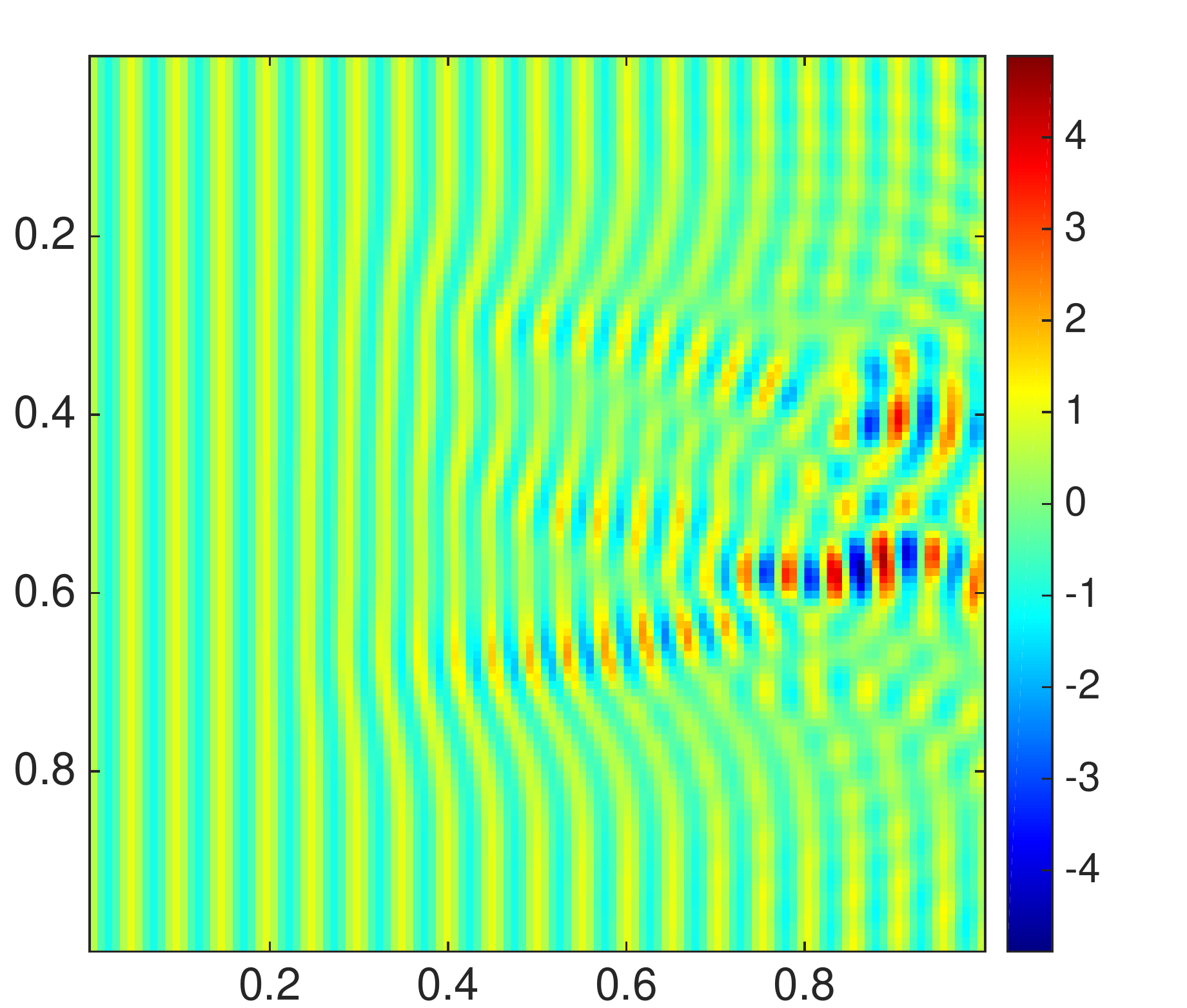}
\begin{tabular}{cc | cc | cc}
\hline
\hline
$k/(2\pi)$ & $N$ & $T_{\text{setup}}$ & $T_{\text{apply}}$ & $N_{\text{iter}}$ & $T_{\text{solve}}$ \\
\hline
$5$ & $3 \times 29^3$ & 2.18e+01 & 6.41e-01 & 6 & 5.62e+00 \\ 
$10$ & $3 \times 59^3$ & 5.08e+02 & 7.21e+00 & 6 & 4.56e+01 \\ 
$20$ & $3 \times 119^3$ & 1.20e+04 & 7.95e+01 & 6 & 5.23e+02 \\ 
\hline
\hline
\end{tabular}
\caption{Results for example \eqref{m:3D3} in 3D. Top left: the inhomogeneity $m(x)$ in
  cross section view at $x_3 = 0.5$. Top right: the third component of the total field $E^i(x) +
  E^s(x)$ at $x_3 = 0.5$ for $k/(2\pi) = 20$. Bottom: the numerical results.}
\label{tab:3D3}
\end{table}

\FloatBarrier

From the numerical results we observe that the iteration number changes at most mildly as the wave
number grows. On the other hand, the iteration number can depend significantly on the profile of
$m(x)$. In the examples, the square/cube obstacle needs more iterations compared to the other
cases. The reason is that the square/cube obstacle has larger areas with high refractive
index. From \eqref{eqn:nbr_i_block_alpha} one can see that larger values of $|m(x)|$ lead to larger
truncation errors and the numerical results are consistent with this observation. Nonetheless, in
all test cases, the iteration numbers are below ten, which show the validity of this preconditioner.

\begin{figure}
[!ht]
\centering
\includegraphics[width=0.45\textwidth]{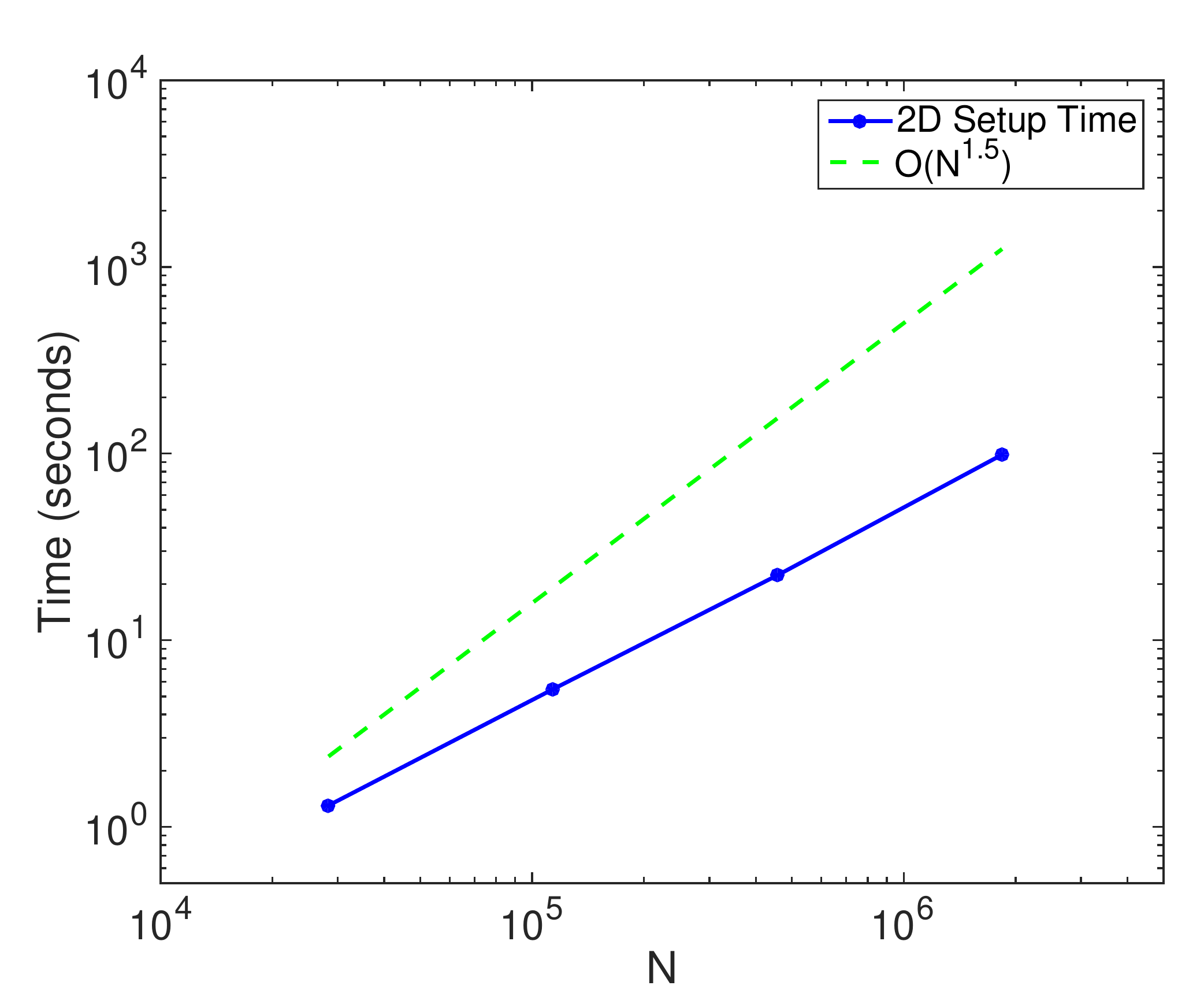}
\quad
\includegraphics[width=0.45\textwidth]{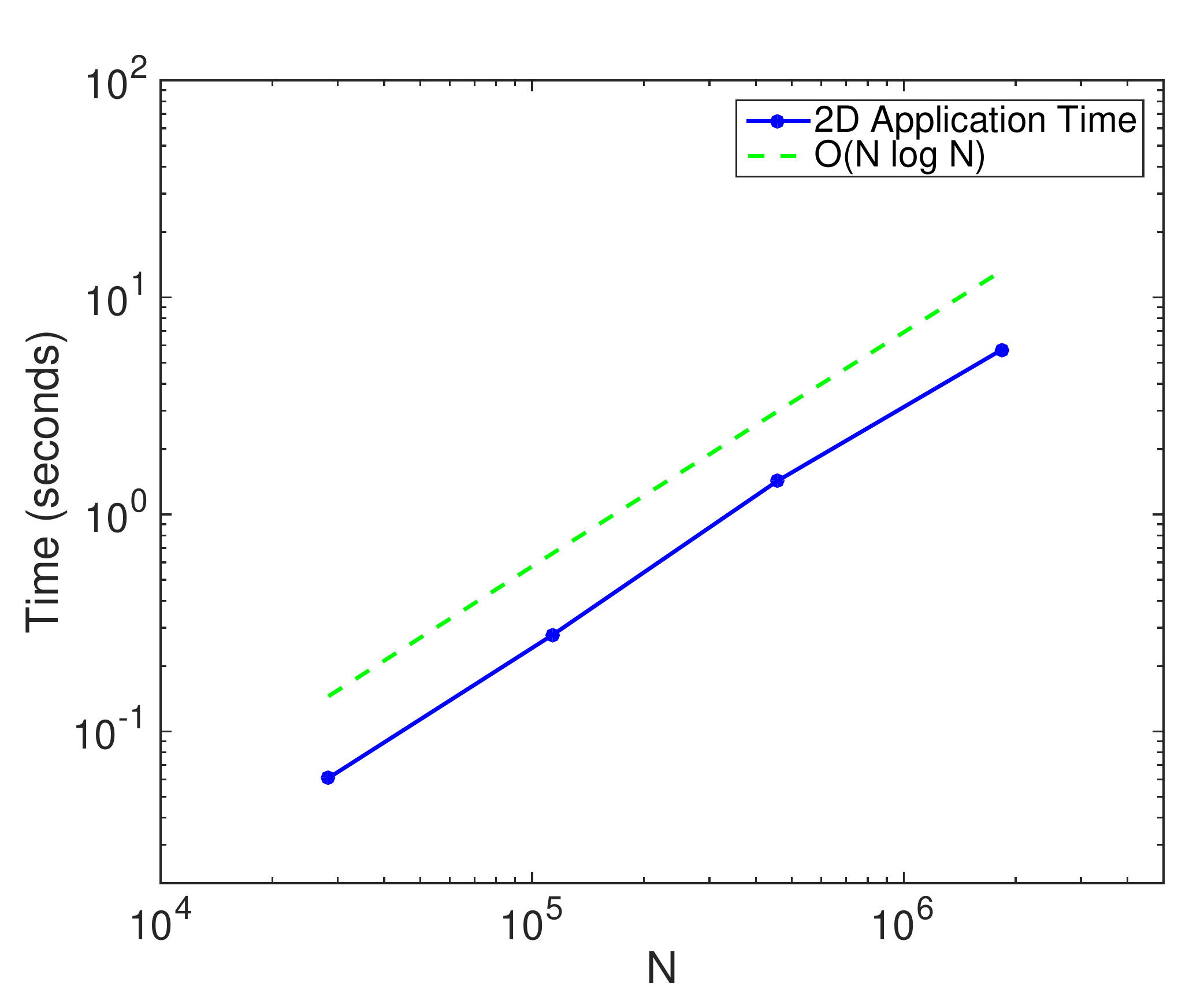}
\caption{The $\log$-$\log$ plots of the scalings of the runtimes in 2D. The runtimes are taken as the averages of $T_{\text{setup}}$ and $T_{\text{apply}}$ of the three test cases. Left: Setup time scaling. Right: Application time scaling. The solid lines are the actual runtimes and the dashed lines are the theoretical scalings. We see that the application time scales as the theoretical cost, while the setup time scales far below, benefiting from MATLAB's built-in parallelization.}
\label{fig:2D_scaling}
\end{figure}

\begin{figure}
[!ht]
\centering
\includegraphics[width=0.45\textwidth]{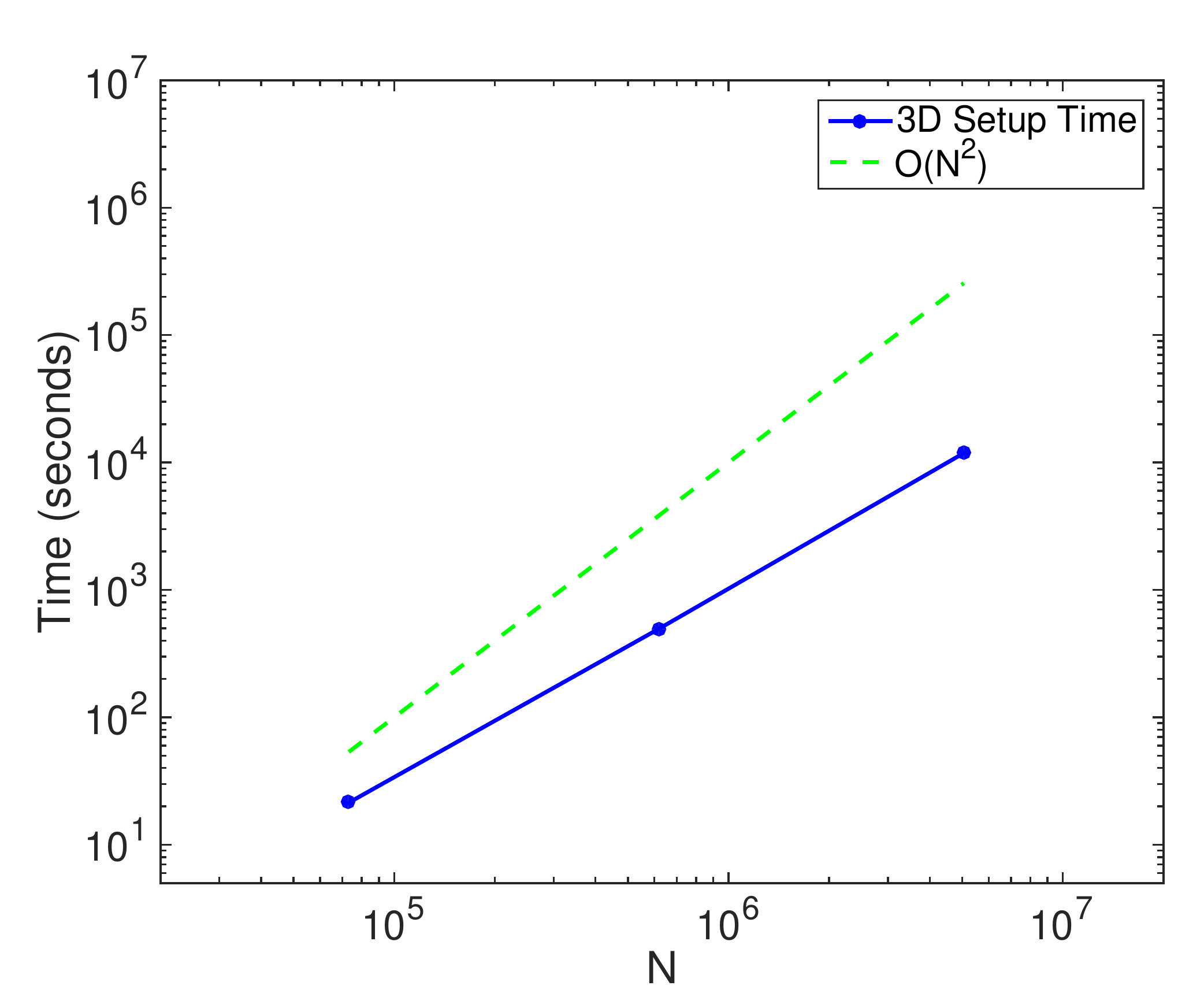}
\quad
\includegraphics[width=0.45\textwidth]{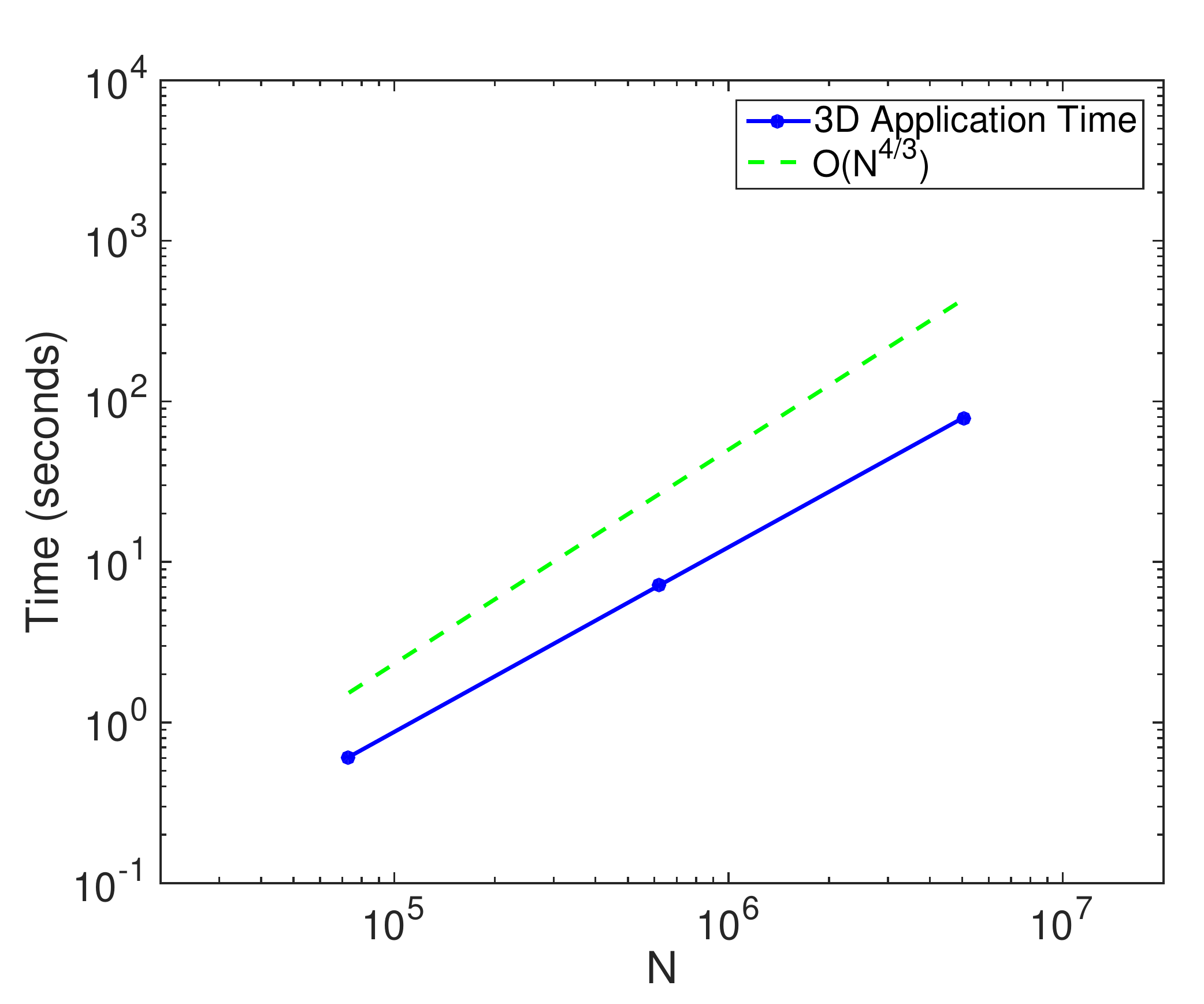}
\caption{The $\log$-$\log$ plots of the scalings of the runtimes in 3D. The runtimes are taken as the averages of $T_{\text{setup}}$ and $T_{\text{apply}}$ of the three test cases. Left: Setup time scaling. Right: Application time scaling. The solid lines are the actual runtimes and the dashed lines are the theoretical scalings. The setup cost scales below the theoretical cost, which is credited to MATLAB's built-in parallelization.}
\label{fig:3D_scaling}
\end{figure}

On the runtime side, the setup and application times are scaling as or below the theoretical complexities, especially in the setup cases where the actual costs are scaling far below the theoretical ones. The credit is to MATLAB's built-in parallelization which notably speeds up the matrix operations. To be specific, it drastically sped up the matrix inversions for the degree of freedoms on the solving front during the setup stage. Figures \ref{fig:2D_scaling} and \ref{fig:3D_scaling} provide the $\log$-$\log$ plot views for setup and application costs in 2D and 3D respectively.

Note that the actual runtime depends on the implementation and platform. If implementing a single thread version, one should hope the costs align closer to the theoretical complexities. What we showed here are two points. First, the runtimes of our implementation scale at least as well as the theoretical analysis. Second, this method can be easily sped up by mature software packages. Ideally, one can use state-of-the-art multifrontal solvers to achieve the best performance.

\section{Conclusions and future work}
\label{sec:conclusion}

This paper presents the sparsifying preconditioner for the time-harmonic Maxwell's equations. The
key idea is to transform the dense linear system into a sparse one by minimizing the non-local
interactions. As shown by the numerical results, when combined with the standard GMRES solver, the
preconditioner converges in only a few iterations, essentially independent of the problem size. The
setup and application costs are almost the same as the ones for solving the sparse system arising
from the PDE formulation.

There are several potential improvements that can be made. First, the problem considered in this
paper only involves inhomogeneity for the electric permittivity $\varepsilon$. Indeed, the magnetic
permeability $\mu$ can be inhomogeneous as well. The only difference is that, by taking the
inhomogeneity of $\mu$ into account, one needs to deal with a larger integral system with both $E$
and $H$ involved. Nevertheless, the same idea applies and the integral system can be sparsified with
a similar procedure. Second, instead of solving the sparsified system with the nested dissection
method, the sweeping preconditioner \cite{ying2011sweephmf,ying2011sweeppml} could be applied to
further reduce the computational cost. There are several previous works indicating the validity of
this approach. For example, \cite{tsuji2012sweepingemfem,tsuji2012sweepingem} applied the sweeping
preconditioner to the time-harmonic Maxwell's equations. \cite{liu2018sparsify} combined the
sweeping preconditioner with the sparsifying preconditioner and formed an efficient preconditioner
that solves the Lippmann-Schwinger equation in quasi-linear time. By combining the two
preconditioners alongside with a recursive approach similar to \cite{liu2016recursive}, we could
hopefully reduce the cost of solving the integral form of the time-harmonic Maxwell's equations to
quasi-linear as well.

\section*{Acknowledgments}
The authors are partially supported by the National Science
Foundation under award DMS-1521830 and the U.S. Department of Energy's
Advanced Scientific Computing Research program under award
DE-FC02-13ER26134/DE-SC0009409.

\bibliographystyle{abbrv}
\bibliography{references}
\end{document}